\documentclass[]{aspm}

\articleinfo{}{}{}

\usepackage{verbatim}
\usepackage{amssymb}
\usepackage{amsbsy}
\usepackage{amscd}
\usepackage{amsmath}
\usepackage{amsthm}
\usepackage[mathscr]{eucal}

\input{mymacros.sty}
\usepackage{pst-all}

\newcommand{\lemref}[1]{Lem\-ma~\ref{#1}}

\def\Eu{\mathop{\mathrm{Eu}}\nolimits}

\newcommand{\Dec}{\operatornamewithlimits{Dec}}

\newenvironment{NB}{
\color{blue}{\bf NB}.  
}{}
\excludeversion{NB}

\title[Wall-crossing formula for framed quiver moduli]{Wall-crossing formula for framed quiver moduli}

\author[R. Ohkawa]{Ryo Ohkawa}

\address{Osaka Metropolitan University, 
3-3-138 Sumiyoshi-ku, Osaka, 558-8585, Japan: ohkawa.ryo@omu.ac.jp, 
Research Institute for Mathematical Sciences, Kyoto University, 
Oiwake-cho, Sakyo-ku, Kyoto, 606-8502, Japan}

\email{ohkawa.ryo@gmail.com}

\rcvdate{}
\rvsdate{}

\subjclass[2010]{}


\begin{document}

\begin{abstract}
We investigate the wall-crossing phenomena for moduli of framed quiver representations. 
These spaces are expected to be highly useful in capturing the representation theoretic essence of special functions 
in integrable systems. 
Within this class of moduli spaces, we focus on the type $A$ flag manifold, type $A$ affine Laumon spaces, 
Nakajima quiver variety, and framed moduli of sheaves on the projective plane and the blow-up
as main motivating examples. 
Specifically, we examine the wall-crossing formulas for integrals of Euler classes over these moduli spaces.
\end{abstract}

\maketitle

\section{Introduction}
\label{sec:intro}
Framed quiver representation is a quiver representation with a framing.
Moduli of framed quiver representations are studied by Reineke \cite{R}. 
These moduli spaces contain interesting examples with applications to various areas such as 
representation theory and integrable systems.
As one of such examples, from affine Laumon space,  
Shiraishi \cite{Shi} introduced non-stationary Ruijsenaars functions.
One of our main results Theorem \ref{thm:adjoint} has framework to study these functions.
On the other hand, affine Laumon spaces also give another type of generating series of
integrals. 
In \cite{AHKOSSY}, via the Shakirov's equation \cite{Shakirov:2021krl} we clarified relationships between $qq$-Painlev\'e equation 
introduced by Hasegawa
\cite{Hasegawa}, \cite{HasegawaLax} and these functions.
As an approach to study these functions systematically, we summarize method using wall-crossing 
formula for framed quiver moduli in this paper. 

Mochizuki \cite{M} studied wall-crossing formula for moduli of parabolic sheaves on surfaces.
He studied integrals over these moduli spaces which depend on stability parameters. 
The set of stability parameters are divided into 
connected open sets called chambers, and the boundary walls.
The wall-crossing formula describes differences between
integrals over two moduli spaces lying on neighboring chambers along one wall
in terms of integrals over Hilbert schemes of points on surfaces.
In particular, it is applied to study Donaldson invariants and Seiberg-Witten invariants of algebraic surfaces, 
but his theory is formulated in vast general setting.

In the case of the blow-up $\hat{\PP}^{2}$ along one point of the projective plane,
Nakajima-Yoshioka \cite{NY2} translated it to quiver language from the context of Nekrasov partition functions.
%
Nekrasov's conjecture \cite{Nek} states that these partition functions give 
deformations of the Seiberg-Witten prepotentials for $N=2$ SUSY Yang-Mills theory.
This conjecture is proven in Braverman-Etingof \cite{BE}, Nekrasov-Okounkov \cite{NO} 
and Nakajima-Yoshioka \cite{NY1} independently.
Originally in \cite{NY1}, they used different method from wall-crossing, 
and only studied the case where integrand is equal to $1$.
But in \cite{NY2}, they also studied integrals of 
various cohomology classes using the method developed by Mochizuki. 
Their blow-up formula are re-formulated as an application of wall-crossing formula 
connecting chambers corresponding to $\PP^{2}$ and $\hat{\PP}^{2}$
deducing the Nekrasov conjecture in more general cases.
This further leads to a proof of the Witten's conjecture 
concerning  Donaldson invariants and Siberg-Witten invariants 
for algebraic surfaces by G\"ottche-Nakajima-Yoshioka \cite{GNY}.

In \cite{O1}, we further studied framed moduli on $\PP^{2}$.
It is difficult to find a useful chamber-and-wall structure
in sheaf theoretic context.
But in terms of quiver, 
we can find two chambers along the wall hyperplane 
perpendicular to imaginary roots. 
This gives a functional equation for Nekrasov partition function in the fundamental matter theory. 
We generalized this computation to the minimal resolution of $A_{1}$-singularity in \cite{O2}. 
We expect this gives another approach to $(-2)$-blowup formula in \cite{S}, \cite{BS1} through 
AGT correspondence \cite{AGT}.
Furthermore applying the method to handsaw quiver variety of type $A_{1}$ in \cite{OY}, 
we obtain geometric interpretation of formulas for multiple hypergeometric functions.
Aparting from moduli spaces of sheaves, we also study the flag manifold of type $A$ in \cite{O3}. 

In this paper, we summarize these techniques in general settings 
of framed quiver moduli.
Here {\it framed quiver} is a pair $(Q, \infty)$ of quiver $Q=(Q_{0}, Q_{1}, Q_{2})$ with relations 
and a vertex $\infty \in Q_{0}$ called 
framing vertex. 

Following ideas in \cite{M} and \cite{NY2}, we give a method to 
describe difference among integrals over framed quiver moduli 
with stability parameters in two chambers along one wall.
As a main result, we obtain a  combinatorial description of wall-crossing term in Theorem \ref{thm:wcm}.
These are analogue of integrals over Hilbert schemes of points on surfaces in Mochizuki theory \cite{M}.
When we take Euler classes of certain direct sum of tautological bundles, this leads to decomposition of 
binomial coefficients.
In particular, we get functional equations \cite[(6)]{O1} of Nekrasov functions, formula \cite[Theorem 3.3]{O2} suggested by Ito-Maruyoshi-Okuda \cite{IMO}, 
and rational limits \cite[Theorem 2.4 (a)]{OY} of the Kajihara-Noumi transformation \cite{K}, \cite{KN}.

Furthermore, we apply Theorem \ref{thm:wcm} to integrals of 
Euler classes of tangent bundles over framed moduli,
and obtain Theorem \ref{thm:adjoint}.
For some framed quivers such as Nakajima quiver varieties, the numbers of arrows connected to
framing vertex $\infty$ do not appear in this theorem. 
For a graph of affine $ADE$, this implies that wall crossing terms are independent of the rank of
framed sheaves as a corollay.
In fact, a few computational experiments suggest that these terms vanish.
But in general, we do not obtain such a result.
When we consider handsaw quiver variety of type $A_{1}$, we get rational limit \cite[Theorem 2.4 (b)]{OY} of transformation 
formula of basic multiple hypergeometric functions 
obtained in Langer-Schlosser-Warnaar \cite{LSW} 
and Halln\"as-Langmann-Noumi-Rosengren \cite{HLNR1}, \cite{HLNR2}
from various contexts.
 
%
%
%

The organization of the paper is the following.
In \S 2, we prepare notation for framed quivers.
In \S 3, we introduce enhancements of quivers and analyze stability conditions of their representations.
In \S 4, we introduce enhanced master spaces of framed quiver representations.
In \S 5, we compute virtual fundamental cycles and normal bundles of fixed points set of enhanced master spaces.
In \S 6, we give a proof for wall-crossing formula in Theorem \ref{thm:wcm}.
In \S 7, we apply Theorem \ref{thm:wcm} to integrals of 
Euler classes of tangent bundles over framed moduli,
and give a proof of Theorem \ref{thm:adjoint}.

The author is grateful for Hidetoshi Awata, Ayumu Hoshino, Akira Ishii, Hiroaki Kanno, 
Hitoshi Konno, Takuro Mochizuki, Kohei Motegi, 
Hiraku Nakajima, 
Masatoshi Noumi, Takuya Okuda,  
Yusuke Ohkubo, Yoshihisa Saito, Jun'ichi Shiraishi, Yuji Terashima, Shintaro Yanagida 
and Yutaka Yoshida for discussion.
He particularly appreciate Takuro Mochizuki for suggesting and encouraging him to write the paper 
in this direction.
He is partially supported by Grant-in-Aid for Scientific Research 21K03180 and 17H06127, JSPS.
He has had the generous support and encouragement of Masa-Hiko Saito.
This work was partly supported by Osaka Central Advanced Mathematical
Institute: MEXT Joint Usage/Research Center on Mathematics and
Theoretical Physics JPMXP0619217849, and by the Research Institute for Mathematical Sciences,
an International Joint Usage/Research Center located in Kyoto University.

\section{Framed Quiver}
\label{sec:quiver}

\subsection{Setting}

\indent We consider a pair $(Q, \infty)$ of a quiver $Q=(Q_{0}, Q_{1}, Q_{2})$ with relations
and $\infty \in Q_{0}$ called {\it a framing vertex}. 
Here $Q_{0}$ is the set of vertices, 
$Q_{1}$ is the set of arrows, and $Q_{2}$ consists of linear combinations of paths with the same beginnings and endings.
Such a pair $(Q, \infty)$ is called {\it framed quiver}, and written as $Q=(Q, \infty)$. 
We put $I = Q_{0} \setminus \lbrace \infty \rbrace$.

For each path $p$, we write by $\text{out}(p)$ and $\text{in}(p)$ the beginning and the ending.
For $l \in Q_{2}$, we also write by $\out(l)$ and $\inn(l)$ the beginning and the ending of paths appearing
in $l$ with non-zero coefficients.
We assume that $\out(l), \inn(l) \in I$ for any $l \in Q_{2}$.

\indent We consider a finite dimensional $Q_{0}$-graded vector space $V=\bigoplus_{v \in Q_{0}} V_{v}$
with $\dim V_{\infty} \le 1$, 
and the set
\[
\rep_{Q}(V)= \prod_{a \in Q_{1}} 
= 
\lbrace
(B_{h}) \in 
\Hom_{\C} (V_{\text{out} (h)}, V_{\text{in} (h)}) \mid
(B_{h}) \text{ satisfies relations in }Q_{2}
\rbrace
\] 
 of $Q$-representations on $V$.
 This is a zero set of a {\it moment map} defined as follows.
 We put
 \[
 \mb M_{Q}(V)=\prod_{a \in Q_{1}} \Hom_{\C} (V_{\out (a)}, V_{\inn (a)}), \quad 
\mb L_{Q}(V)=\prod_{l \in Q_{2}} \Hom_{\C} (V_{\out (l)}, V_{\inn (l)}), 
\]
and define $\mu=\mu_{Q} \colon \mb M_{Q}(V) \to \mb L_{Q}(V)$ by sending
$\rho \in \mb M_{Q}(V)$ to $\mu(\rho)=(\rho(l))_{l \in Q_{2}}$ where $\rho(l) \in \Hom_{\C} (V_{\out (l)}, V_{\inn (l)})$ is defined by 
composing linear maps of $\rho=(B_{a})_{a \in Q_{1}}$ taking linear combinations according to $l$.
 
We consider $\prod_{v \in Q_{0}} \GL(V_{v})$-action on $\rep_{Q}(V)$.
For stability conditions, we take $\zeta = (\zeta_{i} )_{i \in I} \in \mb R^{I}$, and put
\begin{align}
\label{zetainfty}
\zeta_{\infty} =
\begin{cases}
- \sum_{i \in I} \zeta_{i} \dim V_{i} & \text{ when }\dim V_{\infty} = 1\\
0 & \text{ when }\dim V_{\infty} = 0.
\end{cases}
\end{align}
For a $Q_{0}$-graded subspace $S=\bigoplus_{v \in Q_{0}} S_{v}$ of $V$, we put 
$\zeta (S)=\sum_{v \in Q_{0}} \zeta_{v} \cdot \dim S_{i}$.
We alway assume $\zeta(V) =0$.
This automatically holds when $\dim V_{\infty}=1$ by \eqref{zetainfty}. 

\begin{defn}
We say that a $Q$-representation $\rho$ on $V$ is $\zeta$-semistable if for any sub-representation $S$ of $\rho$, we have $\zeta(S) \le 0$.
Furthermore $\rho$ is said to be stable if the inequality is always strict for any non-trivial
proper sub-representation $S$.
\end{defn}

We put $M_{Q}^{\zeta}(V) = \lbrace \rho \in \rep_{Q}(V) \mid 
\zeta\text{-semistable} \rbrace / (\prod_{v \in Q_{0}} \GL(V_{v})/\C^{\ast} \id_{V} )$.
Here we note that $\C^{\ast} \id_{V}$ trivially acts on $\rep_{Q}(V)$.

When $V_{\infty} = \C$, we put 
$G= \prod_{i \in I} \GL(V_{i})$, which can be identified with $(\prod_{v \in Q_{0}} \GL(V_{v})/\C^{\ast} \id_{V} )$ by 
normalizing components in $\GL(V_{\infty})$ to $\id_{V_{\infty}}$.
In the following, we fix $Q$ and mainly consider the case where $V_{\infty}=\C$.
For $\alpha = (\dim V_{v})_{v \in Q_{0}} \in (\Z_{\ge 0})^{Q_{0}}$, we also put
$M^{\zeta}(\alpha)=M_{Q}^{\zeta}(V)$ abbreviating $Q$.

\subsection{Example}
We see some examples of $Q$ with a framing vertex $\infty$, and torus actions acting on $M_{Q}^{\zeta}(V)$.
First we consider the Nakajima quiver variety introduced in \cite{N1}.
\begin{exam}
\label{exam:quivervar}
Let us consider a graph $(I, E)$, where $I$ is the set of vertices, and $E$ is the set of
edges.
We choose an orientation for each edge in $E$, and write by $H$ the set of edges with the orientations.
For each arrow $h \in H$, we write by $\bar{h}$ the reversed arrow, and put $\Omega=H \sqcup \Bar{H}$, where
$\Bar{H}=\lbrace \bar{h} \mid h \in H \rbrace$.
We also write by $\bar{}$ the involution of $\Omega$ sending $h \in H$ to $\bar{h} \in \Bar{H}$. 

For each vector $r=(r_{i})_{i \in I} \in (\Z_{\ge 0})^{I}$, we define a quiver with relations 
$Q=Q^{r}=(Q_{0}, Q_{1}, Q_{2})$
as follows.
We put $Q_{0}=\lbrace \infty \rbrace \sqcup I$, and 
$Q_{1} = \Omega \sqcup
\lbrace 
z_{i\ell}, w_{i \ell} \mid i \in I, \ \ell=1, \ldots, r_{i}
\rbrace$, where $z_{i \ell}$ is an arrow from $\infty$ to $i$, and $w_{i \ell}$ is 
an arrow from $i$ to $\infty$.
For each arrow $a \in \Omega$, we put
\[
\varepsilon(a)
=
\begin{cases}
1 & a \in H \\
-1 & a \in \Bar{H}.
\end{cases}
\] 
Then we put $Q_{2}=\lbrace l_{i} \mid i \in I \rbrace$ where
\[ 
l_{i} = \displaystyle \sum_{\substack{a \in \Omega \\ \out(a) = i }} \e(a) \bar{a} a + \sum_{\ell=1}^{r_{i}} z_{i \ell} w_{i \ell}.
\] 

To define a torus action on $\Rep_{Q}(V)$, we introduce $I$-graded vector space $W=\bigoplus_{i \in I} W_{i}$ with $\dim W_{i} = r_{i}$.
We consider a torus $\mb T = (\C^{\ast} )^{\times 2} \times (\C^{\ast} )^{\times \sum_{i \in I} r_{i}}$, where the second factor is
regarded as a set of diagonal matrices in $\GL(W)$ for certain fixed basis of $W$.
Then we rewrite $\mb M_{Q}(V)$ and $\mb L_{Q}(V)$ as
\begin{align*}
\mb M_{Q}(V)
&=
\bigoplus_{a_{1} \in H} \Hom_{\C}(V_{\out(a_{1})}, V_{\inn(a_{1})}) \otimes \C_{q_{1}} \oplus
\bigoplus_{a_{2} \in \Bar{H}} \Hom_{\C}(V_{\out(a_{2})}, V_{\inn(a_{2})}) \otimes \C_{q_{2}} \\
&\oplus \bigoplus_{i \in I} 
\left\lbrace
\Hom_{\C}(W_{i}, V_{i})
\oplus \Hom_{\C}(W_{i}, V_{i}) 
\otimes \C_{q_{1}q_{2}}
\right\rbrace,\\
\mb L_{Q}(V)
&=
\bigoplus_{i \in I} \End_{\C}(V_{i}) \otimes \C_{q_{1} q_{2}}. 
\end{align*}
These are same as the original vector spaces, but 
can be regarded as $\mb T$-representations naturally.
Since the moment map $\mu \colon \mb M_{Q}(V) \to \mb L_{Q}(V)$ is $\mb T$-equivariant and compatible 
with the $G$-actions, we have the induced $\mb T$-action on $M_{Q}^{\zeta}(V)$.
\end{exam}
We use the following diagram:
\begin{center}
\includegraphics[scale=0.1]{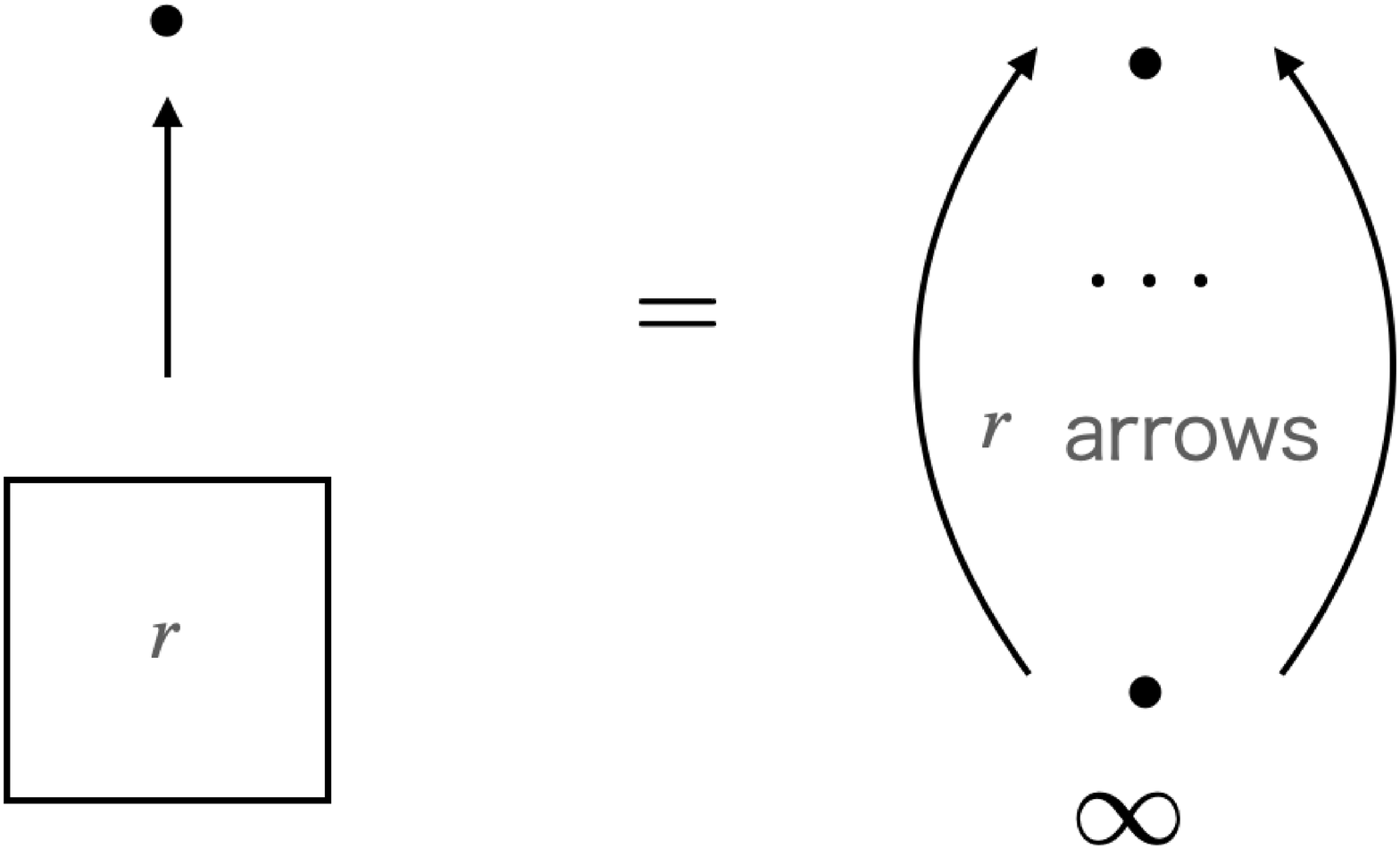}
\hspace{3cm}
\includegraphics[scale=0.1]{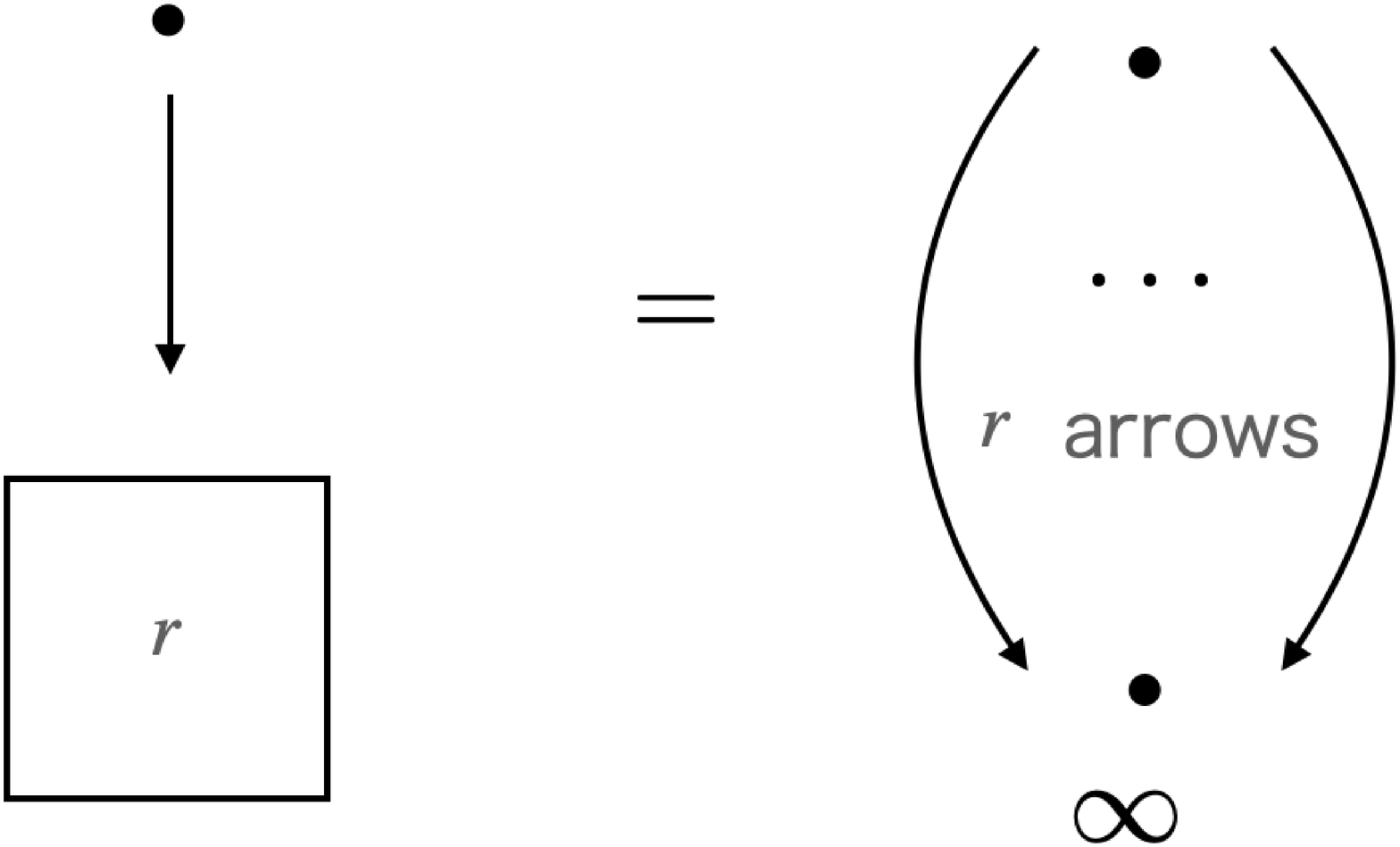}
\end{center}
For example, the quiver variety of type $A_{n}$ is displayed in the following diagram: 
\begin{center}
\includegraphics[scale=0.1]{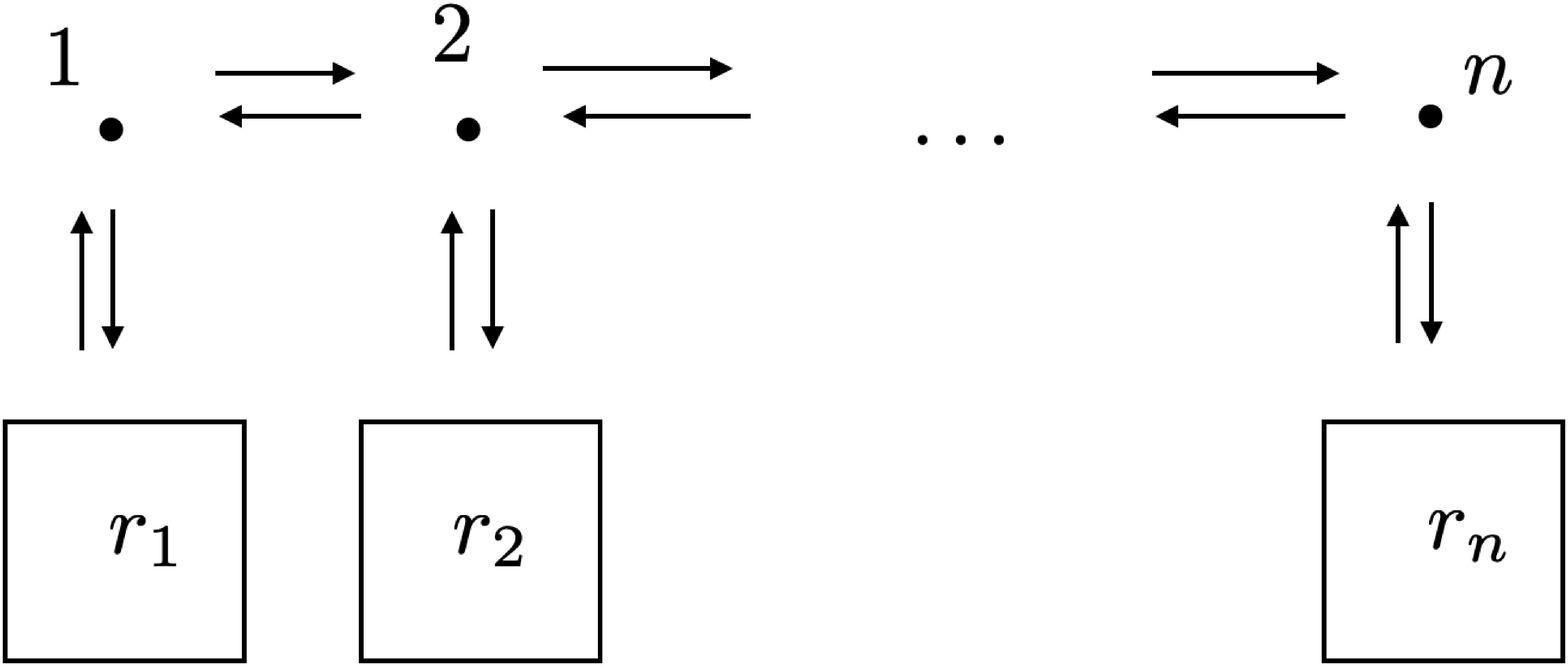}
\end{center}
When $r_{2}=\cdots = r_{n}=0$ and $\zeta_{1}, \ldots, \zeta_{n} < 0$, 
the quiver variety $M_{Q}^{\zeta}(V)$ is isomorphic to the cotangent bundle of 
flag manifold of flags $F_{\bullet}$ in $W_{1}$ with $\dim F^{i} = \dim V_{i}$ for $i=1, \ldots, n$.

\begin{exam}
(Chainsaw quiver variety)
As a vriant of quiver varieties we have chainsaw quiver varieties displayed in the following
diagram:
\begin{center}
\includegraphics[scale=0.1]{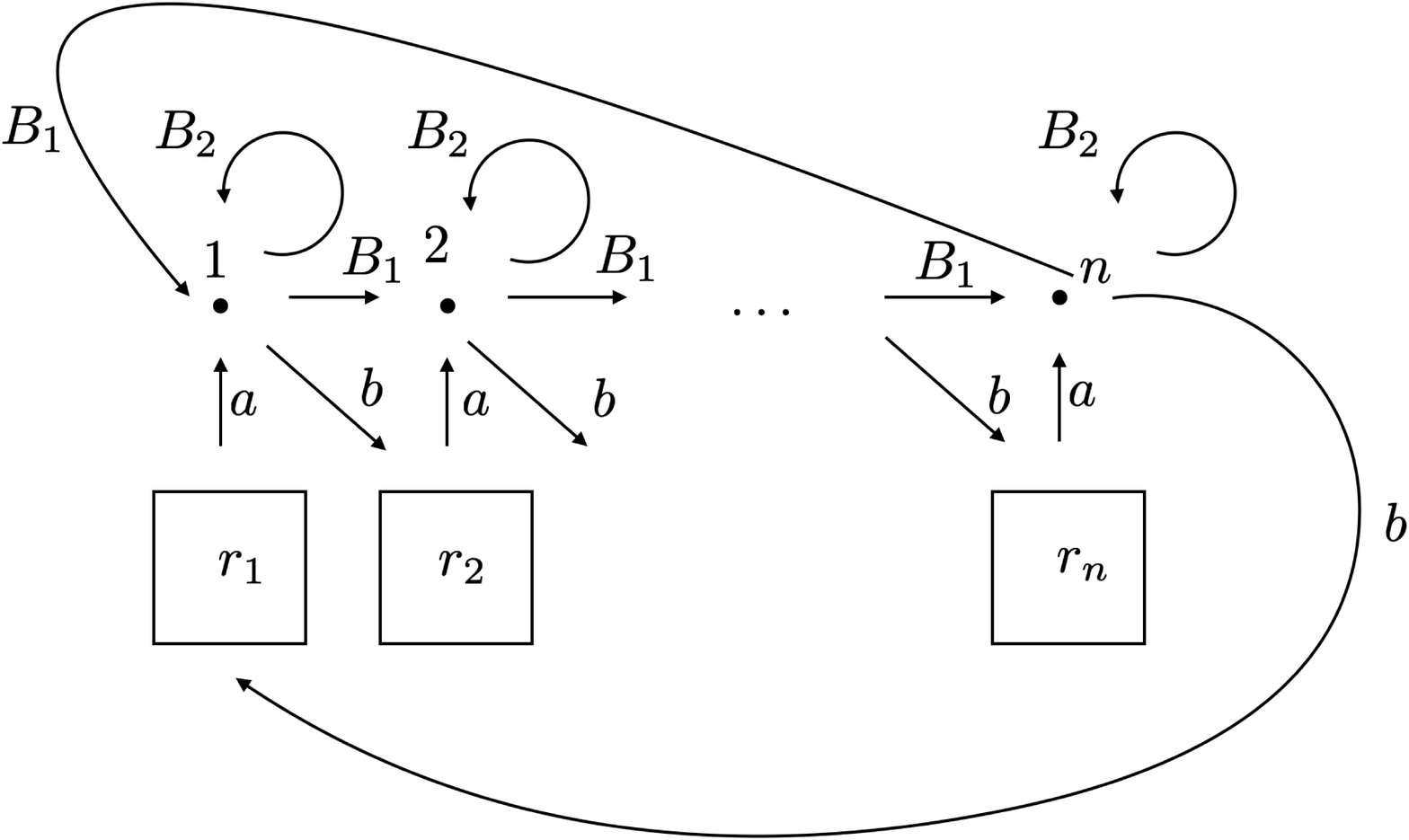}
\end{center}
Here we give relations $[B_{1}, B_{2}] + ab = 0$.
Then we rewrite $\mb M_{Q}(V)$ and $\mb L_{Q}(V)$ as
\begin{align*}
\mb M_{Q}(V)
&=
\prod_{i \in I} \Hom_{\C}(V_{i}, V_{i+1}) \otimes \C_{q_{1}} 
\times
\prod_{i \in I} \Hom_{\C}(V_{i}, V_{i}) \otimes \C_{q_{2}} \\
&\times 
\prod_{i \in I} \Hom_{\C}(W_{i}, V_{i})
\times \Hom_{\C}(W_{i}, V_{i+1}) \otimes \C_{q_{1}q_{2}},\\
\mb L_{Q}(V)
&=
\prod_{i \in I} \End_{\C}(V_{i}) \otimes \C_{q_{1} q_{2}},
\end{align*}
where $I=\Z/n\Z$.
We also have the induced $\mb T$-action on $M_{Q}^{\zeta}(V)$.
\end{exam}

\begin{exam}
\label{exam:blup}
Framed moduli space on the blowup $\hat{\PP}^{2}$ is also obtained from the following diagram:
\begin{center}
\includegraphics[scale=0.1]{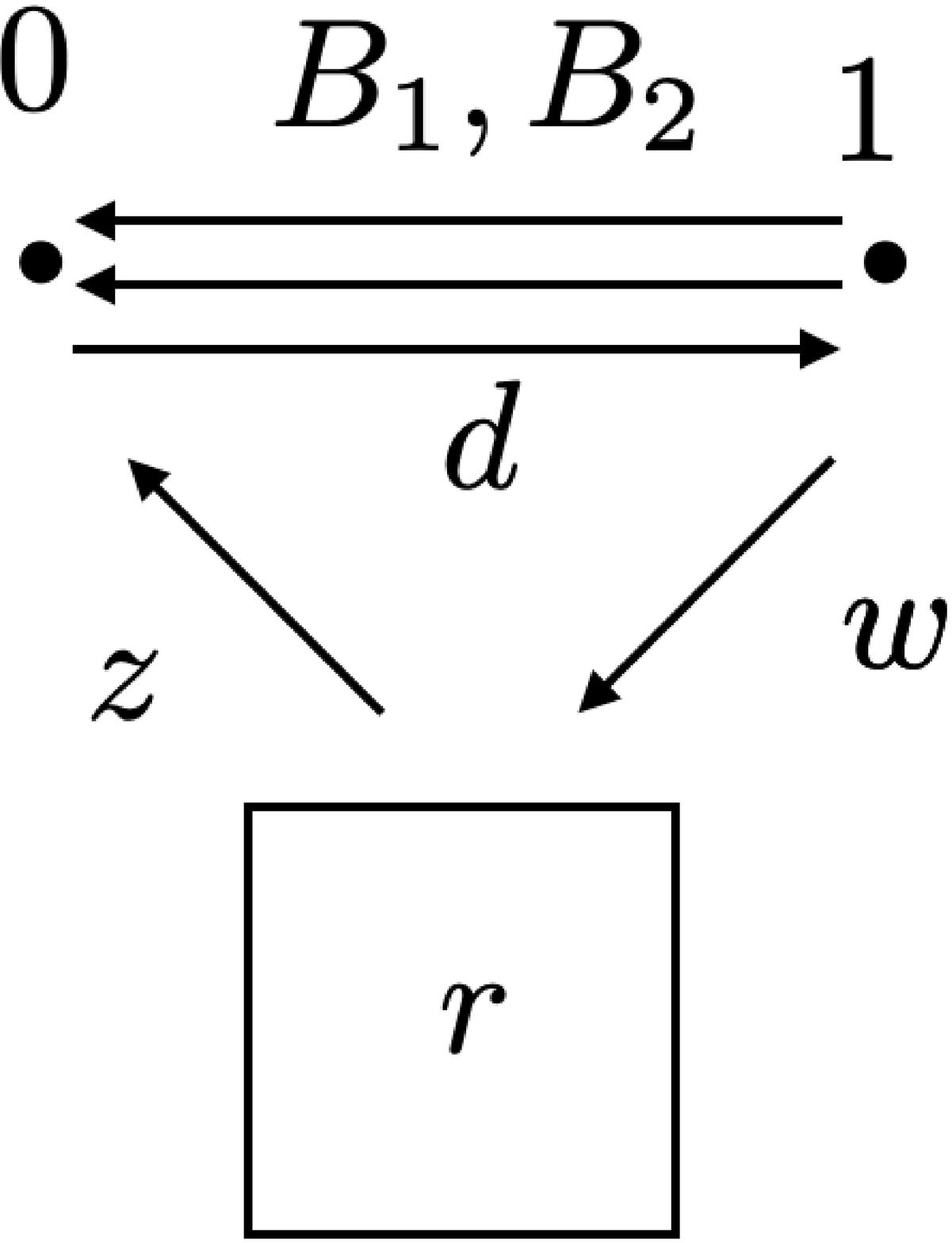}
\end{center}
with relations $B_{1} d B_{2} - B_{2} d B_{1} +zw=0$.
Then we rewrite $\mb M_{Q}(V)$ and $\mb L_{Q}(V)$ as
\begin{align*}
\mb M_{Q}(V)
&=
\Hom_{\C}(V_{0}, V_{1}) \otimes \C_{q_{1}} 
\times
\Hom_{\C}(V_{0}, V_{1}) \otimes \C_{q_{2}} \\
&\times 
\Hom_{\C}(W, V_{0})
\times \Hom_{\C}(V_{1}, W) \otimes \C_{q_{1}q_{2}},\\
\mb L_{Q}(V)
&=
\End_{\C}(V_{0}) \otimes \C_{q_{1} q_{2}} \times \End_{\C}(V_{1}) \otimes \C_{q_{1} q_{2}}.
\end{align*}
We also have the induced $\mb T$-action on $M_{Q}^{\zeta}(V)$.
\end{exam}

\begin{exam}
(Flag manifold)
We also have flag manifolds as important examples displayed in the following diagram:
\begin{center}
\includegraphics[scale=0.1]{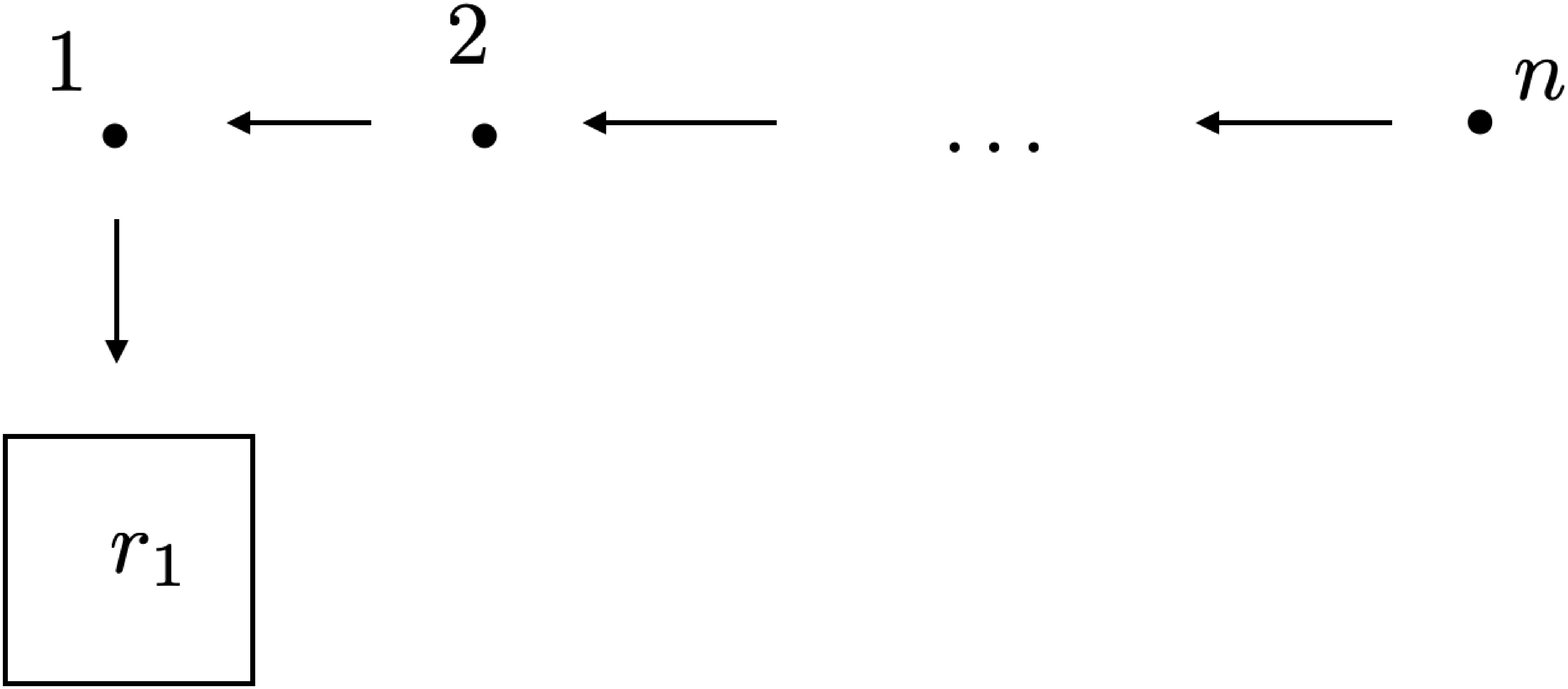}
\end{center}
with no relation, that is, $Q_{2}=\emptyset$.
Then we rewrite $\mb M_{Q}(V)$ as
\begin{align*}
\mb M_{Q}(V)
&=
\prod_{i \in I} \Hom_{\C}(V_{i}, V_{i+1})  
\times 
\Hom_{\C}( V_{1}, W).
\end{align*}
When $\zeta_{1}, \ldots, \zeta_{n} < 0$, we $M_{Q}^{\zeta}(V)$ is isomorphic to the flag manifold.
We also have the induced $\GL(W)$-action on $M_{Q}^{\zeta}(V)$.
They are sub-varieties of quiver varieties of type $A_{n}$. 
\end{exam}

\section{Enhancement of quiver and stability condition}
\label{sec:enhancement}
In the following, we take a $Q_{0}$-graded vector space $V=\bigoplus_{v \in Q_{0}} V_{v}$ such that
$\dim V_{\infty} =1$, and put $\alpha = (\dim V_{i})_{i \in Q_{0}} \in (\mb Z_{\ge 0 })^{Q_{0}}$.

We fix a dimension vector $\beta =(\beta_{i} )_{i \in I} \in (\mb Z_{\ge 0 })^{I}$ such that $\beta_{i} \le \alpha_{i}$
for $i \in I$.
We take $\bar{\zeta} \in \beta^{\perp}=\lbrace \zeta \in \mb R^{I} \mid \sum_{i \in I} \zeta_{i} \beta_{i} = 0 \rbrace$ 
which does not lie on any other hyperplane $(\beta')^{\perp} \neq \beta^{\perp}$ with $\beta_{i}' \le \alpha_{i}$ for
$i \in I$.  
We have two chambers whose boundaries containing $\bar{\zeta}$, and write by $\mc C$ 
the one whose element $\zeta$ satisfy $(\zeta, \beta)<0$, and by $\mc C'$ the other one.
To describe wall-crossing among $\mc C$ and $\mc C'$, we choose $0 \in I$ such that $\beta_{0} \neq 0$, 
and consider enhancement $\tilde{Q}$ of $Q$ as follows.

We analyze behavior of moduli space of $\tilde{Q}$-representations by taking refinement of stability 
parameters $\zeta \in \mc C \sqcup \mc C'$ and $\bar{\zeta} \in \beta^{\perp}$.

\subsection{Enhancement of quiver}
\label{subsec:enhancement}
We choose one vertex $0 \in I=Q \setminus \lbrace \infty \rbrace$ with $\beta_{0} \neq 0$, 
and fix a non-negative integer $L \ge \alpha_{0}$. 
We introduce the following enhancement $\tilde{Q}=(\tilde{Q}_{0}, \tilde{Q}_{1}, \tilde{Q}_{2})$ of our quiver $Q$.
We put $\tilde{Q}_{2}=Q_{2}$.
The set of vertices in $\tilde{Q}$ is 
\[
\tilde{Q}_{0} = Q_{0} \sqcup \lbrace  (0, k) \mid k= 1, 2, \ldots, L \rbrace.
\]
The set $\tilde{Q}_{1}$ of arrows in $\tilde{Q}$ is a disjoint union of $Q_{1}$ and
\[ 
\lbrace (0, k) \to (0, k+1) \mid k = 1, 2, \ldots, L - 1\rbrace \sqcup
\lbrace (0, L) \to 0 \rbrace. 
\]
We also write by $\tilde{Q}$ a framed quiver $(\tilde{Q}, \infty)$.

We consider a $\tilde{Q}_{0}$-graded vector space $\tilde{V}= V \oplus \bigoplus_{k=1}^{L} 
\tilde{V}_{(0,k)}$ with $\dim \tilde{V}_{(0,1)} \le 1$, $\dim \tilde{V}_{(0, L)} = \alpha_{0}$, and
\begin{align}
\label{fullassump}
\dim \tilde{V}_{(0,k-1)} \le \dim \tilde{V}_{(0,k)} \le \dim \tilde{V}_{(0,k-1)} +1 
\end{align}
for $k=2, \ldots, L - 1$.
We put $\mk I=\lbrace k \in [L] \mid \dim \tilde{V}_{(0,k)} - \dim \tilde{V}_{(0,k-1)}= 1\rbrace$.

For a vector bundle $\mc E$ over $M$ and a finite set $\mk I \subset \Z$ with $\rk \E$ elements, 
we write by $Fl_{M}(\mc E, \mk I)$ the full flag bundle 
$F_{\bullet}=(F_{i})_{i \in \Z}$ of $\mc E$ such that 
\[
\lbrace i \in \Z \mid F_{i}/F_{i-1} \neq 0 \rbrace = \mk I.
\]
These are all isomorphic to $Fl_{M}(\mc E, [\rk \mc E])$, but we use products of these flag bundles
and combinatorial description later.
We always put $\mc F_{0}=0$.

\subsection{$(\zeta, \eta)$-stability}

As stability parameters, we take $\zeta \in \mb R^{I}$ and 
$\eta = (\eta_{1}, \ldots, \eta_{L}) \in (\mb Q_{>0})^{L}$, and put
\[
\mu_{(\zeta, \eta)}(\tilde{S})= \frac{ \sum_{v \in Q_{0}} \zeta_{v} \dim S_{i} + \sum_{k=1}^{L} \eta_{k} \dim 
\tilde{S}_{(0,k)} }{\dim S_{\infty} + \sum_{i \in I} \dim S_{i}}
\]
for any $\tilde{Q}_{0}$-graded vector space $\tilde{S}= S \oplus \bigoplus_{k=0}^{L} \tilde{S}_{(0,k)}$
where $S=\bigoplus_{v \in Q_{0}} S_{v}$.
We say that $\tilde{\rho} \in \rep_{\tilde{Q}}(\tilde{V})$ is 
$(\zeta, \eta)$-semistable if for any non-zero proper sub-representation $\tilde{S}$, we have
\begin{eqnarray}
\label{muineq}
\mu_{(\zeta, \eta)}(\tilde{S}) \le \mu_{(\zeta, \eta)}(\tilde{V}).
\end{eqnarray}
If inequality is always strict, we say that $\tilde{\rho}$ is $(\zeta, \eta)$-stable.

\begin{lem}
If $\tilde{\rho}$ is $(\zeta, \eta)$-semistable, then the map $\tilde{V}_{(0, L)} \to V_{0}$ 
is an isomorphism, and 
\[
\tilde{V}_{(0, k)} \to \tilde{V}_{(0, k+1)}  
\] 
are all injective for $k=1, \ldots, L-1$.
\end{lem}
\proof
For any sub-representation $\tilde{S}$ of $\tilde{\rho}$, the inequality \eqref{muineq} is equivalent to
\begin{align}
\label{theta}
\begin{vmatrix}
\sum_{i\in Q_{0}} \zeta_{i} \dim S_{i} + \sum_{k=1}^{L} \eta_{k} \dim \tilde{S}_{(0,k)} &
\sum_{k=1}^{L} \eta_{k} \dim \tilde{V}_{(0,k)}\\
\sum_{v \in Q_{0}} \dim S_{v} &  \dim V
\end{vmatrix}
\le 0.
\end{align}
So we can write the above determinant by 
$\sum_{i \in \tilde{Q}_{0}} \theta_{i} \cdot \dim \tilde{S}_{i}$
with $\theta = (\theta_{i})_{i \in \tilde{Q}_{0}} \in \mb Q^{\tilde{Q}_{0}}$ so 
that $\tilde{\rho}$ is $(\zeta, \eta)$-semistable if and only if 
$\theta(\tilde{S}) \le 0$ for any sub-representation $\tilde{S}$.
Since we have $\theta_{(0,k)} = \eta_{k} \dim V > 0$, 
we get the assertion.
\endproof
We identify $(\zeta, \eta)$-semistable representation $\tilde{\rho}$ with $(\rho, F_{\bullet})$ a $Q$-representation $\rho$ on $V$ with a full flag $F_{\bullet} $ in $V_{0}$
defined by $F_{k}= \im ( \tilde{V}_{(0,k)} \to V_{0})$. 
This flag may have repetitions, but
$\dim (F_{k}/F_{k-1}) \le 1$ by the assumption \eqref{fullassump}.
Under this identification, we have $\mk I=\lbrace k \in [L] \mid \dim F_{k} / F_{k-1}= 1\rbrace$.

To study stability conditions for $(\rho, F_{\bullet})$, it is enough to consider pairs $(S, S_{0} \cap F_{\bullet})$ and $(V/S, F_{\bullet}/ S_{0})$ for all sub-representations $S \subset V$.
Here we write by $F_{\bullet} / S_{0}$ the induced flag $S_{0} + F_{\bullet} / S_{0}$ of $V_{0} / S_{0}$.  
We identify $(S, S_{0} \cap F_{\bullet})$ with the corresponding $\tilde{Q}_{0}$ graded vector space 
$\tilde{S}$.

\subsection{$(\bar{\zeta}, \ell)$-stability condition}

We take a dimension vector
$\beta =(\beta_{i} )_{i \in I} \in (\mb Z_{\ge 0 })^{I}$ and
$\bar{\zeta} \in \beta^{\perp}=\lbrace \zeta \in \mb R^{I} \mid \sum_{i \in I} \zeta_{i} \beta_{i} = 0 \rbrace$ 
which does not lie on any hyperplane other than $\beta^{\perp}$.  
We have two chambers whose boundaries containing $\bar{\zeta}$, and write by $\mc C$ 
the one whose element $\zeta$ satisfy $(\zeta, \beta)<0$, and by $\mc C'$ the other one.

To describe wall-crossing among $\mc C$ and $\mc C'$, we choose $0 \in I$ such that $\beta_{0} \neq 0$, 
and consider enhancement $\tilde{Q}$ as in the previous subsection.
We interpret $(\zeta, \eta)$-stability for suitable choice of $\zeta \in \R^{Q_{0}}$ and 
$\eta = (\eta_{k})_{k=1}^{L} \in 
(\R_{>0})^{L}$.
\begin{defn}
\label{ellstab}
For $\ell \ge 0$, a pair $(\rho, F_{\bullet})$ is said to be $(\bar{\zeta}, \ell)$-stable if $\rho$ is 
$\bar{\zeta}$-semistable and any sub-representation $S$ of $\rho$ with $\bar{\zeta}(S)=0$ satisfies the following two conditions:
\begin{enumerate}
\item[\textup{(1)}] If $S_{\infty}=0$ and $S \neq 0$, we have $S_{0} \cap F_{\ell} =0$.
\item[\textup{(2)}] If $S_{\infty}=\C$ and $S \neq V$, we have $F_{\ell} \not\subset S_{0}$. 
\end{enumerate}
\end{defn}
When $F_{\ell}=0$ (resp. $F_{\ell}= V_{0}$), a pair $(\rho, F_{\bullet})$ is $(\bar{\zeta}, \ell)$-stable if and only if $\rho$ is $\zeta$-stable for $\zeta \in \mc C$ (resp. $\zeta \in \mc C'$).

To construct moduli of $(\bar{\zeta}, \ell)$-stable pairs $(\rho, F_{\bullet})$, 
we consider the following condition:
\begin{align}
\label{conda}
\bar{\zeta} (S) < 0 \ ( \text{resp.} > 0) \implies \mu_{(\zeta, \eta)} (S) < \mu_{(\zeta, \eta)} (V) 
\ ( \text{resp.} > \mu_{(\zeta, \eta)} (V)),
\end{align}
where $S= \bigoplus_{i \in I} S_{i}$ is $I$-graded vector space, but also regarded as $\tilde{Q}_{0}$-graded 
vector space with 
$\tilde{S}_{(0,k)}=0$ for any $k$. 
and the following condition:
\begin{align}
\label{condb}
\frac{\zeta( \beta)}{\sum_{v \in Q_{0}} \beta_{v} } + \sum_{k=\ell+1}^{L} \eta_{k} \cdot \dim F_{k} <
\frac{\sum_{k=1}^{L} \eta_{k} \cdot \dim F_{k} } {\dim V}
< 
\frac{\zeta( \beta)}{\sum_{v \in Q_{0}} \beta_{v} } + \frac{\eta_{\ell} } {\dim V}
\end{align}
for $\ell = 0, 1, \ldots, L$.
Here we put $\eta_{0}=0$.
We remark that $\bar{\zeta}(S)=0$ implies that $(\dim S_{i})_{i \in I}$ is proportional to $\beta \in (\Z_{\ge 0})^{I}$ 
from our generic choice of $\bar{\zeta}$.
 
By the similar arguments as in \cite[\S 4.1]{O1}, we have the following
\begin{prop}
[\protect{
\cite[Proposition 4.2.4]{M}, 
\cite[Lemma 5.6]{NY2}}]
\label{ell}
We take $ \zeta \in \mc C'$, and assume that $( \zeta,  \eta)$ satisfies \eqref{conda} and \eqref{condb}. 
Then for pairs $(\rho, F_{\bullet})$, the $(\bar{\zeta}, \ell)$-stability is equivalent to the $(\zeta, \eta)$-stability.
Furthermore the $(\zeta, \eta)$-semistability automatically implies the $(\zeta, \eta)$-stability.
\end{prop}
\proof
We consider subrepresentation $S$ of $\rho$, and divide into three cases:\\
(i) When $\bar{\zeta}(S)\neq 0$, then from \eqref{conda} inequality $\bar{\zeta}(S) < 0$ is equivalent to 
inequality $\mu_{(\zeta, \eta)}(S, S_{0} \cap F_{\bullet}) < \mu_{(\zeta, \eta)}(V, F_{\bullet})$.\\
(ii) When $\bar{\zeta}(S) =0$ and $S_{\infty} = 0$, we have
\[
\mu_{\zeta, \eta}(S)
\begin{cases}
> \frac{\zeta( \beta)}{\sum_{v \in Q_{0}} \beta_{v}} + \frac{\eta_{\ell}}{\dim V} & ( S_{0} \cap F_{\ell} \neq 0 )\\
< \frac{\zeta( \beta)}{\sum_{v \in Q_{0}} \beta_{v}} + \sum_{k=\ell+1}^{L} \eta_{k} \dim F_{k} 
& ( S_{0} \cap F_{\ell} = 0 ).
\end{cases}
\]
Hence in this case, slope inequality for $(S, S_{0} \cap F_{\bullet})$ is equivalent to $S_{0} \cap F_{\ell} = 0$
by \eqref{condb}.\\
(iii) When $\bar{\zeta}(S) =0$ and $S_{\infty} = \C$, we have
\[
\mu_{\zeta, \eta}(\tilde{V} / \tilde{S})
\begin{cases}
> \frac{\zeta( \beta)}{\sum_{v \in Q_{0}} \beta_{v}} + \frac{\eta_{\ell}}{\dim V} & (  F_{\ell}
\not \subset S_{0} )\\
< \frac{\zeta( \beta)}{\sum_{v \in Q_{0}} \beta_{v}} + \sum_{k=\ell+1}^{L} \eta_{k} \dim F_{k} 
& (  F_{\ell} \subset S_{0} ).
\end{cases}
\]
Hence in this case, slope inequality for $(S, S_{0} \cap F_{\bullet})$ is equivalent to $F_{\ell} \subset S_{0}$
by \eqref{condb}.

In all cases, we get strict inequalities if $(\rho, F_{\bullet})$ is semi-stable.
Hence the last assertion follows.
\endproof

When $F_{\ell} \neq 0$, we can choose $(\zeta, \eta)$ satisfying \eqref{conda}, \eqref{condb} as follows.
First we can take neighborhood of $(\bar{\zeta}, (0, \ldots, 0))$ in $\mb R^{I} \times \mb Q^{L}$ in which any $(\zeta, \eta)$ satisfy \eqref{conda} 
since the coefficients $\theta =( (\theta_{v})_{v \in Q_{0}}, (\theta_{i})_{i \not \in Q_{0}})$ in \eqref{theta} 
approaches to $\dim V \cdot (\bar{\zeta}, 0)$ when $|\zeta - \bar{\zeta}|$ and $\eta$ approches to $0$.
In the following, we search $(\zeta, \eta)$ satisfying \eqref{condb} inside intersection of this neigborhood and $\mc C$. 
We choose $\zeta$ and $\eta_{1}, \ldots, \eta_{\ell} > 0$ such that 
\[
\frac{\sum_{k=1}^{\ell} \eta_{k} \cdot \dim F_{k} - \eta_{\ell} } {\dim V}
< \frac{\zeta( \beta)}{\sum_{i\in Q_{0}} \beta_{i} } < \frac{\sum_{k=1}^{\ell} \eta_{k} \cdot \dim F_{k} } {\dim V}
\] 
Then we take $\eta_{\ell+1}, \ldots, \eta_{L} > 0$ 
small enough so that
\[
\frac{\sum_{k=1}^{L} \eta_{k} \cdot \dim F_{k} - \eta_{\ell} } {\dim V}
< \frac{\zeta( \beta)}{\sum_{i\in Q_{0}} \beta_{i} } < 
\frac{\sum_{k=1}^{\ell} \eta_{k} \cdot \dim F_{k} + \sum_{k=\ell+1}^{L} \eta_{k} \cdot \dim F_{k}(1-\dim V) } {\dim V}
\] 
which is equivalent to \eqref{condb}.

We put $M^{\bar{\zeta},\ell}_{\tilde{Q}}(\tilde{V}) = M_{\tilde{Q}}^{(\zeta, \eta)}(\tilde{V})$ for $\zeta \in \mc C'$ and $\eta \in \mb Q_{>0}^{L}$
satisfying \eqref{conda} and \eqref{condb}, which can be regarded as
moduli of $(\bar{\zeta}, \ell)$-stable pairs $(\rho, F_{\bullet})$ of $Q$-representations $\rho$ on $V$ and full flags $F_{\bullet}$ of $V_{0}$
by the above proposition.
We also put $\wt M^{\ell}(\alpha, \mk I )=M^{\bar{\zeta},\ell}_{\tilde{Q}}(\tilde{V})$.
These are independent of $\mk I \subset [L]$ with $|\mk I|=d$ up to isomorphisms, 
but the index of flags $\mc F_{\bullet}$ depends on $\mk I$.
We also abbreviate $\mk I$ as $\wt M^{\ell}(\alpha )=\wt M^{\ell}(\alpha, \mk I )$ when it is harmless.
 
By remark after Definition \ref{ellstab}, $M^{\bar{\zeta}, \ell}_{\tilde{Q}}(\tilde{V})$ for $\ell =0$ ( resp. $L$ ) 
is the full flag bundle of tautological bundles $\mc V_{0}$ on $M^{\zeta}(\alpha )=M^{\zeta}_{Q}(V)$ for $\zeta \in \mc C$ ( resp. $\zeta \in \mc C'$ ).

\subsection{2-stability conditions}

We consider the following condition on $\eta$:
\begin{equation}
  \label{2stab}
  \sum_{k=1}^{L} \eta_k l_k  \neq 0 \quad
  \text{for any $(l_1,\dots, l_{L} )\in \bigcup_{1 \le m, n \le \alpha_{0} / \beta_{0} } \frac{1}{mn} 
  \Z^{L}\setminus \{0\}$
    with $|l_k|\le \alpha_{0}$}.
\end{equation}
We call \eqref{2stab} $2$-stability condition following \cite{M}.

We remark  again that our flag $F_{\bullet}$ of $V_{0}$ may have repetitions, but assume
$\dim (F_{k}/F_{k-1}) \le 1$, and $F_{0}=0, F_{L}=V_{0}$ as before.

\begin{lem}\label{lem:2-stab}
  Assume \eqref{conda} and \eqref{2stab}.
If $(\rho, F_{\bullet})$ is strictly $(\zeta,\eta)$-semistable, then there
  exists a submodule $0\neq S\subsetneq V$ such that
  \begin{enumerate}
  \item $\mu_{(\zeta,\eta)}(S, S_{0} \cap F_{\bullet}) = \mu_{(\zeta,\eta)}(V, F_{\bullet})$,
  \item $(S, S_{0} \cap F_{\bullet})$ and $(V/S, F_{\bullet} /S_{0})$
    are $(\zeta,\eta)$-stable.
  \end{enumerate}
  Moreover the submodule $S$ is unique except when $(V, F_{\bullet})$ is
  the direct sum 
 \[
  (S, S_{0} \cap F_{\bullet})\oplus (V/S, F_{\bullet} /S_{0}).
  \] 
  In this case the another choice of the submodule is
  $V/S$.
\end{lem}
\proof
  Take a submodule $S$ violating the $(\zeta,\eta)$-stability of
  $(\rho, F_{\bullet})$. Then we have (1). Moreover $(S, S_{0} \cap F_{\bullet})$ and $(V/S,
  F_{\bullet}/ S_{0})$ are $(\zeta,\eta)$-semistable. We have
  either $S_\infty = 0$ or $(V/S)_\infty = 0$.

  Assume either $(S,S_{0} \cap F_{\bullet})$ or $(V/S, F_{\bullet} /S_{0})$ is strictly $(\zeta,\eta)$-semistable. 
  Then we have a
  filtration $0 = V^0 \subsetneq V^1 \subsetneq V^2 \subsetneq V^3 =
  V$ with $\mu_{(\zeta,\eta)}(V^a/V^{a-1},F_{\bullet}^{a}) =
  \mu_{(\zeta,\eta)}(V,F_{\bullet})$ for $a=1,2,3$, where $F_{\bullet}^{a}$
  denote the induced flag $(F_{\bullet} \cap V^a + V^{a-1} )/V^{a-1}$.

  Among $V^a/V^{a-1}$ ($a=1,2,3$), one of them has $\C$ and two of
  them have $0$ at the $\infty$-component. Assume $V^1$ has $\C$ at
  the $\infty$-component for brevity, as the following argument can be
  applied to the remaining cases.

 Then \eqref{conda} implies that $(\dim (V^2/V^1)_{i})_{i \in I}$ and $(\dim (V^3/V^2)_{i} )_{i \in I}$ are 
 proportional to $\beta \in \Z^{I}$.
  Then by the assumption \eqref{2stab},  
the equality  $\mu_{(\zeta,\eta)}(V^2/V^1,F_{\bullet}^{2}) =
  \mu_{(\zeta,\eta)}(V^3/V^2,F_{\bullet}^{3})$
  implies that
$(\dim F_{k}^2)_{k=0}^{L}$ and $(\dim F_{k}^3)_{k=0}^{L}$ are proprotional to each other, 
and hence $(\dim (F_{k}^2/F_{k-1}^{2}) )_{k=0}^{L}$ and $(\dim(F_{k}^3/F_{k-1}^{3}))_{k=0}^{L}$ are so.

On the other hand, we have
\begin{equation*}
  \dim (F_{k}^1/F_{k-1}^1)
  + \dim(F_{k}^2/F_{k-1}^2)
  + \dim(F_{k}^3/F_{k-1}^3) = \dim (F_{k}/F_{k-1}) = 0 \text{ or } 1
\end{equation*}
for any $k$. Therefore at most one of three terms in the left hand
side can be $1$ and the other terms are $0$. 
This implies $\dim(F_{k}^2/F_{k-1}^2) = \dim
(F_{k}^3/F_{k-1}^3) = 0$ for any $k$. 
Thus we get a contradiction $V^2 = V^3$.

If we have another submodule $S'$ of the same property, then
the above argument implies $S\cap S' = 0$ or $S\cap S' = S =
S'$. In the former case we have $S' = V/S$.
\begin{NB}
  Consider $S\cap S'\subset S$ and $\subset S'$, $S/S\cap S'\subset
  V/S'$ and $S'/S\cap S'\subset V/S$.
\end{NB}%
\endproof

\begin{lem}
\label{lem:stabilizer} Let
  $(\zeta,\eta)$ as in \lemref{lem:2-stab}. If
  $(\rho, F_{\bullet})$ is $(\zeta,\eta)$-semistable, its stabilizer is
  either trivial or $\C^*$. In the latter case, $(V,F_{\bullet})$ has
  the unique decomposition $(V^{\flat},F_{\flat}^\bullet)\oplus (V^{\sharp},F_{\sharp}^\bullet)$
  of subrepresentations
  such that both $(V^{\flat},F_{\flat}^\bullet)$, $(V^{\sharp},F_{\sharp}^\bullet)$ are
  $(\zeta,\eta)$-stable, and $\mu_{(\zeta,\eta)}(V^{\flat}, F_{\bullet}^{\flat}) =
  \mu_{(\zeta,\eta)}(V^{\sharp}, F_{\bullet}^{\sharp})=\mu_{(\zeta,\eta)}(V, F_{\bullet})$. The stabilizer comes from that of the
  factor $(V^{\sharp},F_{\sharp}^\bullet)$ with $(V^{\sharp})_\infty = 0$.
\end{lem}
\proof
We take an automorphism $g$ of $(\rho, F_{\bullet})$, and suppose that $g$ has an eigenvalue
  $\lambda\neq 1$. 
  Then we have the generalized eigenspace
  decomposition $(V,F_{\bullet}) = (V^{\flat},F_{\flat}^\bullet)\oplus
  (V^{\sharp},F_{\sharp}^\bullet)$ with $(V^{\flat})_\infty = \C$, $(V^{\sharp})_\infty = 0$.  By
  \lemref{lem:2-stab} $(V^{\flat},F_{\flat}^\bullet)$,
  $(V^{\sharp},F_{\sharp}^\bullet)$ are $(\zeta,\eta)$-stable. Since they have the
  same $\mu_{(\zeta,\eta)}$ and are not isomorphic, there are no nonzero
  homomorphisms between them. Therefore the stabilizer is $\C^*$ in
  this case. The uniqueness follows from that in
  \lemref{lem:2-stab}.

  Next suppose that $g$ is unipotent and let $n= g - 1$. 
  Assume $n\neq 0$
  and let $j \ge 1$ such that $n^j \neq 0$, $n^{j+1} = 0$. 
  We consider
  the submodule $0\neq \ker n^j \subsetneq V$. 
  Since $\ker n^{j}$ and $\im n^{j}$ are both sub-representation of $\rho$, from the
  $(\zeta,\eta)$-semistability of $(\rho, F_{\bullet})$ we have
\(
  \mu_{(\zeta,\eta)}(\ker n^j, \ker n^j \cap F_{\bullet})
  = \mu_{(\zeta,\eta)}(V,F_{\bullet}).
\)
\begin{NB}
From the proof of \lemref{lem:isom}.
\end{NB}
Therefore $(\ker n^j, \ker n^j \cap F_{\bullet})$ and $(V/\ker n^j,
\bar{F}^\bullet/\ker n^j)$ are $(\zeta,\eta)$-stable by
\lemref{lem:2-stab}. They are not isomorphic since they
have different $\infty$-components. 
However $n^{j} \colon V/\ker n^j\to
\ker n^j$ is a nonzero homomorphism, and we have a contradiction.
\endproof

\subsection{Conditions for direct summands}
\label{subsec:Con}
Under the conditions \eqref{conda} and \eqref{2stab}, when $(\zeta, \eta)$-semistable $\tilde{Q}$-representation 
$(\rho, F_{\bullet})$ on $V$ has a non-trivial stabilizer group, we have a decomposition 
$(V, F_{\bullet})= (V^{\sharp}, F_{\bullet}^{\sharp}) \oplus (V^{\flat}, F_{\bullet}^{\flat})$ of 
$\tilde{Q}$-representations by Lemma \ref{lem:stabilizer} with 
$\mu_{(\zeta,\eta)}(V^{\sharp}, F_{\bullet}^{\sharp}) 
= \mu_{(\zeta,\eta)}(V^{\flat}, F_{\bullet}^{\flat})=\mu_{(\zeta,\eta)}(V, F_{\bullet})$. 
Here $V^{\sharp} = \bigoplus_{v \in Q_{0}} V^{\sharp}_{v}$ and 
$V^{\flat} = \bigoplus_{v \in Q_{0}} V^{\flat}_{v}$ are $Q_{0}$-graded vector spaces, and we assume $V^{\flat}_{\infty}=\C$. 
Then we have $(\dim V^{\sharp}_{i})_{i \in I} =p \beta \in \mb Z^{I}$ for an integer $p >0$, 
since the condition \eqref{conda} implies that $\rho$ is $\bar{\zeta}$-semistable
and $\bar{\zeta}(V^{\flat})=\bar{\zeta}(V^{\sharp})$.

We put $\mk I = \lbrace k \in [L] \mid F_{k} / F_{k-1} \neq 0 \rbrace$ as before.
Then the decomposition 
$F_{\bullet} = F_{\bullet}^{\sharp} \oplus F_{\bullet}^{\flat}$
gives us {\it decomposition data} $\mk I^{\sharp} = \lbrace k \in [L] \mid F_{k}^{\sharp} / F_{k-1}^{\sharp} \neq 0 \rbrace
\subset \mk I$ with $|\mk I^{\sharp}|= p \beta_{0}$.

\begin{lem}
\label{lem:minsharp}
Under the conditions \eqref{conda} and \eqref{2stab}, we further assume the first inequality 
in \eqref{condb}.
Then we have $\min (\mk I^{\sharp}) \le \ell$.
\end{lem}
\proof
For the contrary, we assume $\min (\mk I^{\sharp}) > \ell$.
Then by the assumption we have $\mu_{(\zeta, \eta)}(V^{\sharp}, F_{\bullet}^{\sharp} ) < \mu_{(\zeta, \eta)}(V, F_{\bullet} )=
\mu_{(\zeta, \eta)}(V^{\sharp}, F_{\bullet}^{\sharp} )$.
This is a contradiction.
\endproof

In this subsection, we consider another condition for $(\zeta, \eta)$ :
\begin{align}
\label{condc}
\eta_{m} > \dim V \sum_{k=m+1}^{L} \eta_{k} \dim F_{k}  \text{ for } m= 1, 2, \ldots, L. 
\end{align}

\begin{lem}
\label{lem:flatsharp}
Assume \eqref{conda}, \eqref{2stab} and \eqref{condc}.
Then $(\zeta, \eta)$-semistable representation $(\rho, F_{\bullet})$ 
has a non-trivial stabilizer group if and only if
we have a decomposition $(V, F_{\bullet})=(V^{\flat}, F_{\bullet}^{\flat}) \oplus 
(V^{\sharp}, F_{\bullet}^{\sharp})$ of sub-representations of $(\rho,F_{\bullet})$ 
such that the following conditions hold:
\\
(a) $(V^{\flat}, F_{\bullet}^{\flat})$ is $(\bar{\zeta}, \min(I_{\sharp}) -1)$-stable, and \\
(b) $(V^{\sharp}, F_{\bullet}^{\sharp})$ is $\bar{\zeta}$-semistable, and for any sub-representation $S \subset V^{\sharp}$ with $\bar{\zeta}(S)=0$ and $F^{1}_{\sharp} \subset S$, we have $S=V^{\sharp}$. 
\end{lem}
\proof We divide into cases as in the proof of Proposition \ref{ell}.
Since $(\zeta, \eta)$-semistability for $(\rho, F_{\bullet})$ is equivalent to ones for
both $(V^{\flat}, F_{\bullet}^{\flat})$ and  
$(V^{\sharp}, F_{\bullet}^{\sharp})$, it is enough to consider sub-representations $S \subset V^{\flat}$ or
$S \subset V^{\sharp}$.

When $S \subset V^{\flat}$ with $\bar{\zeta}(S) = 0$ and $S_{\infty}=0$, 
$\mu_{(\zeta, \eta)}(S, S_{0} \cap F_{\bullet} ) \le \mu_{(\zeta, \eta)}(V, F_{\bullet})  = \mu_{(\zeta, \eta)}(V^{\sharp}, F_{\bullet}^{\sharp})$ is equivalent to
the following inequality
\[
\frac{\sum_{k=1}^{L} \eta_{k} \cdot \dim S_{0} \cap F_{k}^{\flat}}{\dim S} < \frac{\sum_{k=1}^{L} \eta_{k} \cdot \dim F_{k}^{\sharp}}{\dim V^{\sharp}}.
\]
This is equivalent to $S_{0} \cap F^{\min(I_{\sharp})-1}_{\flat} = 0$ by \eqref{condc}.
The remaining cases follow from the similar arguments.
\endproof

We can take $(\zeta, \eta)$ satisfying \eqref{conda}, \eqref{condb}, \eqref{2stab} and \eqref{condc}
by obvious modifications of the method in the final part of \S 2.3.

\section{Enhanced master space}
\label{sec:enhanced}
In this section, we fix $\ell \in \mk I$. 
We take $\zeta^{-} \in \mc C$, $\zeta^{+} \in \mc C'$ 
and $\eta \in ( \mb Q_{>0} )^{L}$ such that
$(\zeta^{\pm}, \eta)$ satisfies \eqref{conda}, $(\zeta^{+}, \eta)$ satisfies \eqref{condb}, and $\eta$ satisfies \eqref{2stab} and \eqref{condc}.
We put
\[
\zeta^{t} = t \zeta^{+} + (1-t) \zeta^{-}. 
\]
Then $(\zeta^{t}, \eta)$ also satisfies \eqref{conda}, and 
the first inequality in \eqref{condb}. 

\subsection{Enhanced master space}
\label{subsec:enhanced}
We put $G=G_{\tilde{V}}=\prod_{i \in I} \GL(V_{i}) \times \prod_{k=1}^{L} \GL(\tilde{V}_{(0,k)})$, 
and $\wt{\mb M}= \mb M_{\tilde{Q}}(\tilde{V})$.
If necessary we multiply enough divisible positive integer so that we can assume $\zeta^{\pm}$ and $\eta$ are all integers.
We define $\theta_{i}^{\pm} \in \Z$ for $i \in I$ and 
$\theta_{(0,k)}^{\pm} \in \Z$ for $k=1, \ldots, L$ from $\zeta^{\pm}$ and $\eta$ as in \eqref{theta}.
We note $\theta^{+}_{(0,k)}=\theta^{-}_{(0,k)}=\eta_{k} \dim V$.
Then we consider ample $G$-linearizations 
\begin{align}
\label{linearization}
\mc L_{\pm}
&=
\mo_{\mb M} \otimes \C_{\chi_{\theta^{\pm}}}
\end{align}
on $\mb M$.
Here for any character $\chi \colon G \to \C^{\ast}$, we write by $\C_{\chi}$ the weight space, 
that is, one dimensional $G$-representation with the eigenvalue $\chi$, and
for $\theta=( (\theta_{i})_{i \in I}, (\theta_{(0,k)})_{k=1}^{L}) \in \Z^{I} \times \Z^{L}$ we put
$\chi_{\theta} = \prod_{i \in I} (\det_{i})^{\theta_{i}} \times \prod_{k=1}^{L} ( \det_{(0,k)})^{\theta_{(0,k)}}$
for determinants $\det_{i}, \det_{(0,k)} \colon G \to \C^{\ast}$ of corresponding components
$\GL(V_{i}), \GL(\tilde{V}_{(0,k)})$.

We put
$\wh{\mb M} = \wh{\mb M}_{\tilde{Q}}(\tilde{V}) = 
\Proj \text{Sym} (\mc L_{-} \oplus \mc L_{+})$ and consider the semi-stable
locus $\wh{\mb M}^{ss}$ with respect to $\mo_{\wh{\mb M}}(1)$.
We define a enhanced moment map $\hat{\mu} \colon \wh{\mb M} \to \mb L_{Q}(V)$ by composition of 
the projection to $\mb M_{Q}(V)$ and $\mu \colon \mb M_{Q}(V) \to \mb L_{Q}(V)$.
We define an enhanced master space by $\mc M = [\hat{\mu}^{-1}(0)^{ss}/ G]$, where $\hat{\mu}^{-1}(0)^{ss} = \hat{\mu}^{-1}(0) 
\times_{\wh{\mb M}} 
\wh{\mb M}^{ss}$.
We consider $\C^{\ast}_{\hbar}$-action on $\mc M$ defined by
\begin{align}
\label{act0}
(\rho, F_{\bullet},  [x_{-}, x_{+}]) \mapsto (\rho, F_{\bullet}, [e^{\hbar} x_{-}, x_{+}]),
\end{align}
and the fixed point sets $\mc M^{\C^{\ast}_{\hbar}}$.
We have obvious components 
$\mc M_{\pm} = [\hat{\mu}^{-1}(0)^{ss} \times_{\wh{\mb M}} \Proj \text{Sym} (\mc L_{\pm}) / G]$
of $\mc M^{\C^{\ast}_{\hbar}}$.
We see that $\mc M_{-}$ is 
isomorphic to the full flag bundle $Fl(\mc V, \mk I)$ of the tautological bundle $\mc V_{0}$ on $M^{\zeta^{-}}_{Q}(V)$
with $\lbrace k \in [L] \mid \mc F_{k} / \mc F_{k-1} \rbrace = \mk I$, and $\mc M_{+}$ is isomorphic to 
$\wt M^{\ell}(\alpha )=M^{\bar{\zeta},\ell}_{\tilde{Q}}(\tilde{V})$.

We take $(\rho, F_{\bullet}, [x_{0}, x_{1} ] )\in \hat{\mu}^{-1}(0)$ representing
a point in $\mc M^{\C^{\ast}_{\hbar}} \setminus \mc M_{+} \sqcup \mc M_{-}$.
By Lemma \ref{lem:stabilizer}, we have a direct sum doecomposition
$V=V^{\sharp} \oplus V^{\flat}$ of $Q_{0}$-graded vector spaces 
with $(\dim V^{\sharp}_{i} )_{i \in I} = d^{\sharp} \beta \in \Z^{I}$ for $d^{\sharp}>0$, $V^{\sharp}_{\infty}=0$ and 
$V^{\flat}_{\infty}=\C$ compatible with $\rho$ and $F_{\bullet}=F_{\bullet}^{\sharp} \oplus F_{\bullet}^{\flat}$
such that $F_{\bullet}^{\sharp}, F_{\bullet}^{\flat}$ are full flags
of $V_{\sharp}, V^{\flat}$ respectively.
We put $\mk I^{\sharp}=\lbrace k \in [L] \mid \dim F^{\sharp}_{k} / F^{\sharp}_{k-1}= 1\rbrace$.
We have $\mk I^{\sharp} \subset \mk I$, $|\mk I^{\sharp}|= d^{\sharp} \beta_{0}$, and 
$\min(\mk I^{\sharp}) \le \ell$ by Lemma \ref{lem:minsharp}.

Motivated by this observation, we introduce the set 
\begin{align}
\label{decompdata}
\mc D^{\ell}(\mk I) = \lbrace \mk I^{\sharp} \subset \mk I \mid |\mk I^{\sharp}| = d^{\sharp} \beta_{0} \text{ for } 
d^{\sharp} > 0, \text{ and } 
\min(\mk I^{\sharp}) \le \ell \rbrace,
\end{align}
of {\it decomposition data}.
We identify an element $\mk I^{\sharp} \in \mc D^{\ell}(\mk I)$ with a pair $(\mk I^{\sharp}, \mk I^{\flat})$, where 
$ \mk I^{\flat} = \mk I \setminus \mk I^{\sharp}$.

\subsection{Modified $\C^{\ast}_{\hbar}$-action}
\label{subsec:modified}
From a decomposition datum $\mk I^{\sharp} \in \mc D^{\ell}(\mk I)$, 
we describe a connected component $\mc M_{\mk I^{\sharp}}$ of 
$\mc M^{\C^{\ast}_{\hbar}} \setminus \mc M_{+} \sqcup \mc M_{-}$ as follows. 
We take a decomposition $\tilde{V}=\tilde{V}^{\sharp} \oplus \tilde{V}^{\flat}$ of a $\tilde{Q}_{0}$-graded vector space
such that $V \cap \tilde{V}^{\sharp}=V^{\sharp}, V \cap \tilde{V}^{\flat} = V^{\flat}$, and
$\mk I^{\sharp} =\lbrace k \in [L] \mid 
\dim \tilde{V}^{\sharp}_{(0,k)} - \dim \tilde{V}^{\sharp}_{(0,k-1)} = 1\rbrace$.

We put $D=\sum_{i \in I} (\theta^{+}_{i} - \theta^{-}_{i} ) \beta_{i}$, and consider a modified action 
$\C^{\ast}_{e^{\hbar/d^{\sharp}D}} \times \mc M \to \mc M$ induced by 
\begin{align}
\label{act}
\left( \rho, F_{\bullet}, [x_{-},x_{+}] \right) \mapsto 
\left( e^{\hbar/d^{\sharp}D} \id_{\tilde{V}^{\sharp} \oplus } \id_{\tilde{V}^{\flat}} \right) 
\left( \rho, F_{\bullet}, [e^{\hbar}x_{-},x_{+}] \right). 
\end{align}
This action is equal to the original $\C^{\ast}_{\hbar}$-action \eqref{act0}, since the difference is absorbed in $G$-action.
We consider a moment map $\tilde{\mu} = \mu_{\tilde{Q}}$ from $\mb M_{\tilde{Q}}(\tilde{V}^{\sharp})$ and  
$\mb M_{\tilde{Q}}(\tilde{V}^{\flat})$.
We write by $\tilde{\mu}^{-1}(0)^{\sharp}$ subsets of 
$\tilde{\mu}^{-1}(0) \cap \mb M_{\tilde{Q}} (\tilde{V}^{\sharp})$ 
consisting of elements satisfying the conditions in Lemma \ref{lem:flatsharp} (b),
and by $\tilde{\mu}^{-1}(0)^{\flat}$ is the $(\bar{\zeta}, \min(\mk I^{\sharp})-1)$-stable locus of 
$\tilde{\mu}^{-1}(0) \cap \mb M_{\tilde{Q}}(\tilde{V}^{\flat})$.
By Lemma \ref{lem:flatsharp}, we have an open embedding
\begin{align}
\label{decomp0}
\left[
\tilde{\mu}^{-1}(0)^{\sharp} \times \left( \C_{\chi_{\sharp} t^{-1} } \setminus \lbrace 0 \rbrace \right)
\times
\tilde{\mu}^{-1}(0)^{\flat} \times \left( \C_{ \chi_{\flat } t} \setminus \lbrace 0 \rbrace \right) \Big/
G_{\tilde{V}^{\sharp}} \times G_{\tilde{V}^{\flat}} \times \C^{\ast}_{t} 
\right]
\to 
\mc M,
\end{align}
where $\chi_{\sharp}$ and $\chi_{\flat}$ are restrictions of $\chi_{\theta^{+}} / \chi_{\theta^{-} }$ to
$G_{\tilde{V}^{\sharp}}$ and $G_{\tilde{V}^{\flat}}$ respectively.
If we write by $\mc M_{\mk I^{\sharp}}$ its image, then we have $\mc M^{\C^{\ast}_{\hbar} } \setminus
\mc M_{+} \sqcup \mc M_{-}=\bigsqcup_{\mk I^{\sharp} \in \mc D^{\ell}(\mk I)} \mc M_{\mk I^{\sharp}}$.

\indent If we replace the group $\C^{\ast}_{t}$ with $\C^{\ast}_{t^{1/d^{\sharp}D}}$, then we have an \'etale cover 
$\Phi_{\mk I^{\sharp}} \colon \mc M'_{\mk I^{\sharp}} \to \mc M_{\mk I^{\sharp}}$ of degree $1/d^{\sharp}D$. 
Further by replacing $V^{\sharp}$ with $V^{\sharp} \otimes \C_{t^{1/d^{\sharp}D}}$ in \eqref{decomp0},  
the component $\C_{\chi_{\sharp}t}$ changes to
\[
\bigotimes_{i \in I} 
\left(
\det V^{\sharp}_{i} \otimes \C_{t^{1/d^{\sharp}D}}
\right)^{\otimes \theta^{+} - \theta^{-}} \otimes \C_{t}
=\C_{\chi_{\sharp}}
\]
while the other components remain.
Hence we get an isomorphism
$\mc M'_{\mk I^{\sharp}} \cong \mc M_{\sharp} \times \mc M_{\flat}$ with
\begin{align}
\mc M_{\sharp} &= \left[ 
\tilde{\mu}^{-1}(0)^{\sharp} \times \left( \C_{\chi_{\sharp} } \setminus \lbrace 0 \rbrace \right)
\Big/ G_{\tilde{V}^{\sharp}} \right],
\label{msharp}
\\
\mc M_{\flat} &= \left[
\tilde{\mu}^{-1}(0)^{\flat} \times \left( \C_{ \chi_{\flat } t} \setminus \lbrace 0 \rbrace \right) 
\Big/ G_{\tilde{V}^{\flat}} \times \C^{\ast}_{t^{1/d^{\sharp}D}}
\right].
\label{mflat}
\end{align}

\subsection{Moduli stack of destabilizing objects $\mc M_{\sharp}$}
\label{subsec:moduli}
To study $\mc M_{\sharp}$, we give a new quiver $Q^{\sharp}=(Q^{\sharp}_{0}, Q^{\sharp}_{1}, Q_{2}^{\sharp})$ 
from $Q=(Q_{0}, Q_{1}, Q_{2})$ as follows.
We put $Q^{\sharp}_{0}=Q_{0} \sqcup \lbrace \infty' \rbrace=I \sqcup 
\lbrace \infty \rbrace \sqcup \lbrace \infty' \rbrace$, $Q^{\sharp}_{1}= Q_{1} \sqcup \lbrace \infty' \to 0 \rbrace$ 
adding a new arrow $\infty' \to 0$, and
$Q_{2}^{\sharp} = Q_{2}$.
We write by $Q^{\sharp}$ a framed quiver $(Q^{\sharp}, \infty')$.
We only consider $Q^{\sharp}_{0}$-graded vector space 
$V^{\sharp} \oplus V^{\sharp}_{\infty'}$ where $V^{\sharp} =\bigoplus_{i \in I}V^{\sharp}_{i}$
such that $(\dim V^{\sharp}_{i})_{i \in I} = d^{\sharp} \beta \in \Z^{I}$ as in \S \ref{subsec:enhanced} and 
$
V^{\sharp}_{\infty'}=\C.
$ 

For any parameter $\zeta \in \mb R^{I}$, we define 
stability parameter $\zeta^{\sharp}=(\zeta^{\sharp}_{v})_{v \in Q^{\sharp}_{0}} \in \mb R^{Q^{\sharp}_{0}}$ 
by
\begin{align}
\label{zetasharp}
\begin{cases}
\zeta^{\sharp}_{i} = \zeta_{i} \text{ for }i \in I \\
\zeta^{\sharp}_{\infty}=0\\
\zeta^{\sharp}_{\infty'} = - \sum_{i \in I} \zeta_{i} \dim V^{\sharp}_{i}=- \zeta(d^{\sharp} \beta).
\end{cases}
\end{align}
We consider chamber $\mc C \subset \mb R^{I}$ with the boundary containing $\bar{\zeta}$ 
as in \S \ref{sec:enhancement}. 
For $\zeta \in \mc C$, we have $\bar{\zeta}(S) \cdot \zeta(S) \ge 0$ for any $I$-graded 
vector space $S \subset V^{\sharp}$, 
and $\zeta(\beta)<0$.
Hence we have the following lemma.
\begin{lem}
\label{lem:sharp}
For $\zeta^{\sharp}$ defined by \eqref{zetasharp} from $\zeta \in \mc C$, 
a representation $\rho_{\sharp} \in \text{Rep}_{Q^{\sharp}}(V^{\sharp} \oplus V^{\sharp}_{\infty'})$ is 
$\zeta^{\sharp}$-stable if and only if it is
$\bar{\zeta}$-semistable, and for any sub-representation $S \subset V^{\sharp} \oplus V^{\sharp}_{\infty'}$ 
with $\bar{\zeta}(S)=0$ and $S_{\infty'}=\C$, we have $S=V^{\sharp} \oplus V^{\sharp}_{\infty'}$. 
\end{lem}
We put $H_{Q}(d^{\sharp} \beta)=M^{\zeta}_{Q^{\sharp}} (V^{\sharp} \oplus V^{\sharp}_{\infty'})$
for $\zeta \in \mc C$,
and consider a tautological homomorphism $\mc V^{\sharp}_{\infty'} \to \mc V^{\sharp}_{0}$ 
on $H_{Q}(d^{\sharp} \beta)$ corresponding to the arrow from $\infty'$ to $0$.
It is injective from $\zeta^{\sharp}$-stability.
We put $\bar{\mc V}^{\sharp}_{0} =\bar{\mc V}^{\sharp}_{0}$ and 
$\bar{\mk I}^{\sharp} =\mk I^{\sharp} \setminus \lbrace \min(\mk I^{\sharp}) \rbrace$, 
and consider the full flag bundle $Fl_{H_{Q}(d^{\sharp} \beta)}
(\bar{\mc V}^{\sharp}_{0}, \bar{\mk I}^{\sharp} )$ of $\mc F_{\bullet}^{\sharp}=(\mc F_{k}^{\sharp} )$
with $\lbrace k \in [L] \mid \mc F_{k}^{\sharp} / \mc F_{k-1}^{\sharp} \rbrace = \bar{\mk I}^{\sharp}$.
We write by the same letter $\mc V^{\sharp}_{0}$ the pull back of $\mc V^{\sharp}_{0}$ to 
$Fl_{H_{Q}(d^{\sharp} \beta)}(\bar{\mc V}^{\sharp}_{0}, \bar{\mk I}^{\sharp})$, 
and consider the determinant line bundle 
$\bigotimes_{i \in I} (\det \mc V^{\sharp}_{i} )^{\otimes (\theta^{+}_{i} - \theta^{-}_{i})}$ over 
$Fl_{H_{Q}(d^{\sharp} \beta)}(\bar{\mc V}^{\sharp}_{0}, \bar{\mk I}^{\sharp} )$.
Then for $\mc M_{\sharp}$ defined in \eqref{msharp}, we have the following proposition.
\begin{prop}
\label{prop:dest}
Over $H_{Q}(d^{\sharp} \beta)$, we have an isomorphism  
\[
\mc M_{\sharp} \cong 
\left[
\left( 
\bigotimes_{i \in I}  \left( \det \mc V^{\sharp}_{i} \right)^{\otimes (\theta^{+}_{i} - \theta^{-}_{i})} 
\otimes \C_{u^{d^{\sharp}D}} 
\right)^{\times} \Big/ \C^{\ast}_{u}
\right].
\]
In particular, we have an \'etale morphism $\mc M_{\sharp} 
\to Fl_{H_{Q}(d^{\sharp} \beta)}(\bar{\mc V}^{\sharp}_{0}, \bar{\mk I}^{\sharp})$
of degree $1/(d^{\sharp}D)$.
\end{prop}
\proof
In \eqref{msharp}, 
we tensor $V^{\sharp}_{(0, \min(\mk I^{\sharp}))}$ with 
any graded piece of $V^{\sharp}$ other than $V^{\sharp}_{(0, \min(\mk I^{\sharp}))}$ itself.
Then the composition 
$\C=V^{\sharp}_{\infty'} \to V^{\sharp}_{(0,\min(\mk I^{\sharp})+1)} 
\to \cdots \to V^{\sharp}_{(0,N)} \to  V^{\sharp}_{0}$
induces a map 
$ \tilde{\mu}^{-1} (0) \to \Rep_{Q_{\sharp}} (V^{\sharp} \oplus V^{\sharp}_{\infty'})$.
By comparing Lemma \ref{lem:flatsharp} (b) and Lemma \ref{lem:sharp}, the semistable locus
$\tilde{\mu}^{-1}(0)^{\sharp}$ coincide with the pull-back of 
$\Rep_{Q_{\sharp}}^{\zeta^{\sharp}} (V^{\sharp} \oplus V^{\sharp}_{\infty'})$ via this map.
We take a natural group automorphism of $G_{\tilde{V}^{\sharp}}$ such that
the above replacement using tensor product
is $G_{\tilde{V}^{\sharp}}$-equivariant via this automorphism.
Then we have a morphism $\mc M_{\sharp} \to H_{Q}(d^{\sharp} \beta) =
M^{\zeta}_{Q^{\sharp}} (V^{\sharp} \oplus V^{\sharp}_{\infty'})$.
This induce the desired isomorphism since we change $\C_{\chi_{\sharp}}$ to
$\C_{\chi_{\sharp} (\det_{(0, \min(\mk I^{\sharp})})^{d^{\sharp} D}}$ via the replacement. 
\endproof

For tautological bundles $\mc V^{\sharp} \oplus \mc V^{\sharp}_{\infty'}$ and $\mc V^{\flat}$ over 
$H_{Q} (d^{\sharp} \beta )$
and $\wt M^{\min(\mk I^{\sharp})-1}(\alpha - d^{\sharp} \beta, \mk I^{\flat})
$, we also
write by the same letter their pull-backs to 
$Fl_{H_{Q}(d^{\sharp} \beta)}(\bar{\mc V}^{\sharp}_{0}, \bar{\mk I}^{\sharp}) \times 
\wt M^{\min(\mk I^{\sharp})-1}(\alpha - d^{\sharp} \beta, \mk I^{\flat})$
by the projections.
By \eqref{msharp}, \eqref{mflat} and Proposition \ref{prop:dest}, we have a natural \'etale morphism 
$\Psi_{\mk I^{\sharp}} \colon \mc M'_{\mk I^{\sharp}} \to 
Fl_{H_{Q}(d^{\sharp} \beta)}(\bar{\mc V}^{\sharp}_{0}, \bar{\mk I}^{\sharp} ) \times 
\wt M^{\min(\mk I^{\sharp})-1}(\alpha - d^{\sharp} \beta, \mk I^{\flat})$ of degree $1/(d^{\sharp}D)^{2}$.
We take universal full flags $\bar{\mc F}_{\bullet}^{\sharp}$ on 
$Fl_{H_{Q} (d^{\sharp} \beta )}( \bar{\mc V}^{\sharp}_{0}, \bar{\mk I}^{\sharp})$, 

On $\mc M_{\mk I^{\sharp}}$, we have also tautological bundle $\mc V^{\sharp}, \mc V^{\flat}$ 
written by the same letter, and 
$\mc V^{\sharp}_{(0,k)}, \mc V^{\flat}_{(0,k)}$ for $k=1, \ldots, L$ 
via the description \eqref{decomp0},

\subsection{Decomposition of $\mc M^{\C^{\ast}_{\hbar}}$}
Summarizing, we have the following theorem.
\begin{thm}
\label{decomp}
We have
\begin{align}
\label{decomp1}
\mc M^{\C^{\ast}_{\hbar}} = \mc M_{+} \sqcup \mc M_{-} \sqcup 
\bigsqcup_{\mk I^{\sharp} \in \mc D^{\ell}(\mk I) }\mc M_{\mk I^{\sharp}}
\end{align}
such that the followings hold.\\
(i) We have $\mc M_{+} \cong \wt M^{\ell}(\alpha, \mk I)$ and $\mc M_{-}$ is isomorphic to the full flag bundle 
$Fl(\mc V_{0}, \mk I)$ of the tautological bundle $\mc V_{0}$ on $M^{\zeta^{-}}(\alpha)$
with $\lbrace k \in [L] \mid \mc F_{k} / \mc F_{k-1} \rbrace = \mk I$.
\\
(ii) For each $\mk I^{\sharp} \in \mc D^{\ell}(\mk I)$, we have finite \'etale morphisms 
$
\Phi_{\mk I^{\sharp}} \colon \mc M'_{\mk I^{\sharp}} \to \mc M_{\mk I^{\sharp}}$ of degree $1/(d^{\sharp} D)$, and 
$\Psi_{\mk I^{\sharp}} \colon \mc M_{\mk I^{\sharp}}' \to  
Fl_{H_{Q}(d^{\sharp} \beta)}( \bar{\mc V}^{\sharp}_{0}, \bar{\mk I}^{\sharp}) \times 
\wt M^{\min(\mk I^{\sharp})-1}(\alpha  - d^{\sharp} \beta, \mk I^{\flat})
$ 
of degree $1/(d^{\sharp}D)^{2}$, 
where $d^{\sharp}=|\mk I^{\sharp}|/\beta_{0}$, $\mk I^{\flat}=\mk I \setminus \mk I^{\sharp}$, and $D=\sum_{i \in I} (\theta^{+}_{i} - \theta^{-}_{i} ) \beta_{i}$. \\
(iii) As $\C^{\ast}_{e^{\hbar/d^{\sharp}D}}$-equivariant 
$Q_{0}$-graded vector bundles on $\mc M_{\mk I^{\sharp}}'$, 
we have a decomposition
\[
\Phi_{\mk I^{\sharp}}^{\ast} \mc V|_{\mc M_{\mk I^{\sharp}}} 
\cong
\left( \Psi_{\mk I^{\sharp}}^{\ast} \mc V^{\sharp} \otimes L 
\otimes e^{\hbar/d^{\sharp}D} \right)
\oplus
\Psi_{\mk I^{\sharp}}^{\ast} \mc V^{\flat}
\]
where $L$ is a line bundle on $\mc M_{\mk I^{\sharp}}'$ such that 
$
L^{\otimes d^{\sharp}D} \cong 
\Psi_{\mk I^{\sharp}}^{\ast} \bigotimes_{\beta_{i} \neq 0 } 
( \det \mc V^{\flat }_{i} \otimes \det \mc V^{\sharp}_{i})^{\otimes (-\theta^{+}_{i} + \theta^{-}_{i})}$. 
Furthermore we have $\Phi_{\mk I^{\sharp}}^{\ast} \mc V^{\sharp}_{(0, \min(\mk I^{\sharp}))} 
\cong L \otimes e^{\hbar/d^{\sharp}D}$, 
$\Phi_{\mk I^{\sharp}}^{\ast} \left( \mc V^{\sharp}_{(0, k)} /  \mc V^{\sharp}_{(0,k-1)} \right) 
\cong \Psi_{\mk I^{\sharp}}^{\ast} \bar{\mc F}^{\sharp}_{k} \otimes e^{\hbar/d^{\sharp}D}$, 
and $\Phi_{\mk I^{\sharp}}^{\ast} \mc V^{\flat}_{(0, k)} 
\cong \Psi_{\mk I^{\sharp}}^{\ast} \mc V^{\flat}_{(0, k)}$ for $k=1, \ldots, N$.
\end{thm}
\proof
The isomorphisms in (i) and (ii) is obtained in \S \ref{subsec:modified} and \S \ref{subsec:moduli}.
For (iii), we summarize the description \eqref{msharp}, \eqref{mflat} and 
Proposition \ref{prop:dest} as follows.
We replace $\C^{\ast}_{t}$ with $\C^{\ast}_{t^{1/kD}}$ and 
tensor $\C_{t^{1/d^{\sharp}D}} \otimes V^{\sharp}_{(0, \min(\mk I^{\sharp}))}$
with $\tilde{V}^{\sharp} / V^{\sharp}_{(0,\min(\mk I ^{\sharp}))}$ and tensor $\C_{t^{1/d^{\sharp}D}}$ with 
$V^{\sharp}_{(0,\min(\mk I^{\sharp}))}$.
This replacement is $G_{\tilde{V}^{\sharp}}$-equivariant via a natural group automorphism of 
$G_{\tilde{V}^{\sharp}}$.
Hence we get an a degree $k D$ \'etale cover  $\Phi_{\mk I^{\sharp}}$ and
$\Phi_{\mk I^{\sharp}}^{\ast} \tilde{\mc V}^{\sharp}/ \mc V^{\sharp}_{(0, \min(\mk I^{\sharp})} \cong 
\left( \Psi_{\mk I^{\sharp}}^{\ast} \tilde{\mc V}^{\sharp}/ \mc V^{\sharp}_{(0, \min(\mk I^{\sharp})} 
\right) \otimes \C_{t^{1/d^{\sharp} D}} \otimes 
\Phi_{\mk I^{\sharp}}^{\ast} \mc V^{\sharp}_{(0, \min(\mk I^{\sharp}))} 
$. 

On the other hand, $\C_{\chi_{\flat} t^{1/d^{\sharp} D}}$ and $\C_{\chi_{\sharp}}\otimes 
\Phi_{\mk I^{\sharp}}^{\ast} \mc V^{\sharp}_{(0, \min(\mk I^{\sharp}))}$ trivialize
via a degree $k^{2} D^{2}$ \'etale cover $\Psi \colon \mc M'_{k} \to M^{-}(W, V^{\flat}) 
\times M^{-}(V_{(1,1)}, V_{\sharp})$ 
by replacing $\det_{(0, \min(\mk I^{\sharp}))}, t^{1/k}$ with $s=\det_{(0, \min(\mk I^{\sharp}))}^{d^{\sharp}D}, t$.
Hence if we put $L=\C_{t^{1/d^{\sharp} D}} \otimes 
\Phi_{\mk I^{\sharp}}^{\ast} \mc V^{\sharp}_{(0, \min(\mk I^{\sharp}))}$, then 
we have $L^{\otimes d^{\sharp} D} \cong \Psi_{\mk I^{\sharp}}^{\ast} \bigotimes_{\beta_{i} \neq 0 } 
( \det \mc V^{\flat }_{i} \otimes \det \mc V^{\sharp}_{i})^{\otimes (-\theta^{+}_{i} + \theta^{-}_{i})}$. 
This implies the isomorphism of tautological bundles.

For $\C^{\ast}_{\hbar}$-action on $\mc M_{k}$, we normalize $x_{+} /x_{-}=1$ in \eqref{act}.
Then we need a gauge transformation $e^{\hbar/d^{\sharp}D} \id_{\tilde{V}_{\sharp}}$.
\endproof

We write by $N_{+}, N_{-}$, and $N_{\mk I^{\sharp}}$ normal bundles of $\mc M_{+}, \mc M_{-}$ and $\mc M_{\mk I^{\sharp}}$ in $\mc M$ respectively.
In the following, we describe these normal bundles using the obstruction theory.

\section{Integral}
\label{sec:obst}
In this section we define virtual fundamental cycle and apply localization formula.
We compute normal bundle for later calculation.

In this paper, constructions of moduli stacks are essentially fitted in \cite[\S 2]{GP}, 
but slightly different point is that we include group actions as in \cite[\S 2.4]{M}.
We can easily modify the argument since we have intersection theory developed in 
\cite{Kr} similar to one in \cite{F}.


\subsection{Virtual fundamental cycle}
\label{subsec:virt}
We recall constructions of virtual fundamental cycles \cite{BF} in the following setting.
We consider a triple $\mc K=(Y, P, \sigma)$. 
Here $Y= (Y, \mo_{Y}(1))$ is a smooth scheme with an ample line bundle $\mo_{Y}(1)$, $P$ is a group scheme acting on $(Y, \mo_{Y}(1))$, and $\sigma \colon Y \to E$ is a $P$-invariant section of $P$-equivariant vector bundle $E$ over $Y$.
We take a semistable locus $Y^{ss}$ of $Y$ with respect to $\mo_{Y}(1)$, and put $Z=\sigma^{-1}(0) \cap Y^{ss}$.
We assume that any point in $Y^{ss}$ has a finite stabilizer group.
Hence $\mc Y=[Y^{ss}/P]$ is a Deligne-Mumford stack.

We define a moduli stack $\mc Z=\mc Z_{\mc K}=[Z/P]$ as the quotient stack.
We call $\mc K=(Y, P, \sigma)$  
a {\it Kuranishi chart} of $\mc Z$.
In our case, $E$ is always a trivial bundle $Y \times \C^{l}$, and we have $\sigma=\id_{Y} \times \varphi \in \Gamma(Y, E)$
where $\varphi \colon Y \to \C^{l}$ is a $P$-invariant function. 
In such a case, we just write $\mc K=(Y, P, \varphi)$ instead of $\sigma = \id_{Y} \times \varphi$.

To define {\it virtual fundamental cycle} $[\mc Z]^{vir}$, we consider 
a vector bundle $\mc E=[(E|_{Y^{ss}}) / P]$ over $\mc Y$.
We consider the zero section $0 \colon \mc Y \to \mc E$ 
and write by the same letter $\sigma \colon \mc Y \to \mc E$
the section induced by $\sigma \in \Gamma(Y, E)$. 
We also consider a fundamental cycle $[\mc Y]$ of smooth Deligne-Mumford stack $\mc Y$.
For example, we can take a Kuranishi chart $\mc K_{\mc M}=(\wh{\mb M}, G, \hat{\mu})$ 
of the enhanced master space $\mc M$ 
as we saw in the previous section.

Using these we define the virtual fundamental cycle of $\mc Z$ by $[\mc Z]^{vir}=0^! [\mc Y]$, 
where $0^!$ is the refined Gysin homomorphism by the zero section $0 \colon \mc N \to \mc E$
in the following diagram :
\[
\xymatrix{
\mc Y \ar[r]^{\sigma} & \E \\
\ar[u] \mc Z \ar[r] & \mc Y \ar[u]^{0}},
\]

We take another Kuranishi chart $\mc K_{1}=(Y_{1}, P_{1}, \sigma_{1})$ of another stack 
$\mc Z_{1}=\mc Z_{\mc K_{1}}$ where $\sigma_{1} \in \Gamma(Y_{1}, E_{1})$ and $E_{1}$
is a vector bundle on $Y_{1}$. 
We assume that there exists a bundle morphism 
$E_{1} \to E$ equivariant via a group homomorphism 
$P_{1} \to P$ such that sections $\sigma_{1}$ and $\sigma$ coincide 
via the isomorphism $E_{1} \to E|_{Y_{1}}$.
This collection of morphisms are called a morphism of Kuranishi charts.
Furthermore when these morphisms induces an \'{e}tale morphism $[Y_{1}^{ss}/P_{1}] \to [Y^{ss}/P]$, 
we call such a collection of morphisms a {\it coordinate change} of Kuranishi charts. 

When we have a morphism $\mc K_{1}=(Y_{1}, P_{1}, \sigma_{1}) \to 
\mc K=(Y, P, \sigma)$ of Kuranishi charts, such that $Y_{1} \to Y$ is flat, then we have 
$[\mc Z_{1}] = \Psi^{\ast} [\mc Z]$.
In general, virtual fundamental cycles can be different depending on choice of Kuranishi chatrts. 


We remark that we can always take Kuranishi chart globally in our setting.
In general, we need obstruction theory developed in \cite{BF} to construct virtual fundamental cycle.

\subsection{Localization formula}
\label{subsec:local1}
In Kuranishi chart $\mc K=(Y, P, \sigma)$, we suppose that $Y$ and $E$ admit $( \C^{\ast} ) ^{\ell}$-actions
compatible with $P$-actions such that the section $\sigma$ is equivariant.
Then $( \C^{\ast} ) ^{\ell}$ acts on $\mc Z$ and $[\mc Z]^{vir}$ defines a class
of $( \C^{\ast} ) ^{\ell}$-equivariant cohomology group $A_{\bullet}^{( \C^{\ast} ) ^{\ell}} (\mc Z)$.
We take a decomposition $\mc Y^{( \C^{\ast} ) ^{\ell}} = \bigsqcup_{\alpha=1}^{m} \mc Y_{\alpha}$
of connected components.

By the localization formula \cite[Corollary 5.3.6]{Kr} for smooth Deligne Mumford stacks, we have
\begin{align}
\label{localsmooth}
[\mc Y] = \sum_{\alpha=1}^{m} \frac{j_{\alpha \ast} [\mc Y_{\alpha}]}{ \Eu(N_{\mc Y_{\alpha}})}
\end{align}
where $j_{\alpha} \colon \mc Y_{\alpha} \to \mc Y$ is the embedding and $N_{\mc Y_{\alpha}}$
is the normal bundle of $\mc Y_{\alpha}$ in $\mc Y$.
\begin{NB}
For the diagram
\[
\xymatrix{
\mc Y_{\alpha} \ar[r]^{j_{\alpha}} & \mc Y \ar[r]^{e} & \E \\
\ar[u] \mc Z_{\alpha} \ar[r]^{\iota_{\alpha}} &\ar[u] \mc Z \ar[r] & \mc Y \ar[u]^{0}},
\]
we have $0^{!} j_{\alpha \ast} = \iota_{\alpha \ast} 0^{!}$ by \cite[Theorem 6.2 (a)]{F}, 
where $0 \colon \mc Y \to \E$ is the zero-section.
Furthermore for the diagram
\[
\xymatrix{
\mc Y_{\alpha} \ar[r]^{j_{\alpha}} & \E_{\alpha} \ar[r] & \E_{\alpha} \oplus (\E|_{\mc Y_{\alpha}})^{mov} \\
\ar[u] \mc Z_{\alpha} \ar[r]^{\iota_{\alpha}} &\ar[u]^{0_{\alpha}} \mc Y_{\alpha} \ar@{=}[r] & 
\mc Y_{\alpha} \ar[u]^{0|_{\mc Y_{\alpha}}}},
\]
we have $0^{!}= e((\E|_{\mc Z_{\alpha}})^{mov}) \cap 0_{\alpha}^{!}$ 
by the excess intersection formula \cite[Theorem 6.3]{F}, 
where $0_{\alpha} \colon \mc Y_{\alpha} \to \E_{\alpha}$ is the zero-section.
\end{NB}
We consider $\mc E|_{\mc Y_{\alpha}}$ with fiber-wise $\mb T$-action, and 
a decomposition $\mc E|_{\mc Y_{\alpha}}= \mc E_{\alpha} \oplus (\mc E|_{\mc Y_{\alpha}})^{mov}$
where $\mc E_{\alpha}$ is the sub-bundle with the trivial weight and 
$(\mc E|_{\mc Y_{\alpha}})^{mov}$ is the sub-bundle with non-trivial weights.  
We put $\mc Z_{\alpha} = \mc Z \times_{\mc Y} \mc Y_{\alpha}$ and define the virtual fundamental
cycle by $[\mc Z_{\alpha}]^{vir}= (0_{\alpha})^{!} [\mc Y_{\alpha}]$ 
where $0_{\alpha} \colon \mc Y_{\alpha} \to \E_{\alpha}$ is the zero-section.

Applying the refined Gysin homomorphism $0^{!}$ to \eqref{localsmooth}, we get
\begin{align}
\label{localDM}
[\mc Z]^{vir} = \sum \frac{\iota_{\alpha \ast} [\mc Z_{\alpha}]^{vir}}{\Eu (N_{\mc Z_{\alpha}})}
\end{align}
in $A_{\bullet}^{( \C^{\ast} ) ^{\ell}} (\mc Z)$.
Here $N_{\mc Z_{\alpha}}=
N_{\mc Y_{\alpha}}|_{\mc Z_{\alpha}} - (\E|_{\mc Z_{\alpha}})^{mov}$ is called 
{\it normal bundle} of $\mc Z_{\alpha}$ in $\mc Z$, and 
$(\E|_{\mc Z_{\alpha}})^{mov}$ is the moving part
of $\E|_{\mc Z_{\alpha}}$.
We note that $N_{\mc Z_{\alpha}}$ depends on Kuranishi charts $\mc K$.

\subsection{Components in $\mc M^{\C^{\ast}_{\hbar}}$}
\label{subsec:comp}

In the following, we consider a Kuranishi chart $\mc K_{\mc M}=\left( \wh{\mb M}, G, \hat{\mu} \right)$ 
of the enhanced master space $\mc M$. 
We have $\C^{\ast}_{\hbar}$-action defined in \eqref{decomp2} which is compatible with $G$-action.
Furthermore we assume that the Kuranishi chart $\mc K$ admits 
$\mb{T}=(\C^{\ast})^{\ell}$-action compatible with 
$\C^{\ast}_{\hbar} \times G$-action such that $\hat{\mu}$ is also 
$\mb{T} \times \C^{\ast}_{\hbar}$-invariant. 

In this setting, we can get Kuranishi charts for $\mc M^{\C^{\ast}_{\hbar}}$ explicitly.  
In fact, the arguments so far hold for
objects in the ambient space $\mb M$, $\wt{\mb M}$ and $\wh{\mb M}$, 
since we can forget $Q_{2}$ from $Q=(Q_{0}, Q_{1}, Q_{2})$.
The ambient stack $\mc N= [\wh{\mb M}^{ss}/G]$ of the enhanced master space 
$\mc M=[\hat{\mu}^{-1}(0) /G]$ is a smooth Deligne-Mumford stack with $\mb{T} \times \C^{\ast}_{\hbar}$-action compatible with one on $\mc M$.
By Theorem \ref{decomp} we have a decomposition
\begin{align}
\label{decomp2}
\mc N^{\C^{\ast}_{\hbar}}= \mc N_{+} \sqcup \mc N_{-} \sqcup 
\bigsqcup_{\mk I^{\sharp} \in \mc D^{\ell}(\mk I)} \mc N_{\mk I^{\sharp}},
\end{align}
where $\mc N_{\pm}=[\Proj \text{Sym} (\mc L_{\pm}) / G]$, and $\mc N_{\mk I^{\sharp}}$ 
is the closed substack of $\mc N^{\C^{\ast}_{\hbar}}$ of elements 
over direct sum $(X_{\flat}, F_{\bullet}^{\flat}) \oplus (X_{\sharp}, F_{\bullet}^{\sharp})$ 
with the decomposition type $\mk I^{\sharp} \in \mc D^{\ell}(\mk I)$.
\begin{NB}
For a scheme $S$, the full sub-category $\mc N_{\mk I^{\sharp}}(S)$ consists of objects $\phi \in \mc N(S)$ such that the pull-back of vector fields corresponding to $\C^{\ast}_{\hbar}$-action vanish, and for any closed point $p \in S$, the corresponding data $(X, F_{\bullet}, [\rho, 1])$ satisfies that $(X, F_{\bullet})$ is $(X_{\flat}, F_{\bullet}^{\flat}) \oplus (X_{\sharp}, F_{\bullet}^{\sharp})$ with the decomposition type $\mk I$.
\end{NB}

In the following, for $\gamma=\pm$, or $\mk I^{\sharp} \in \mc D^{\ell}(\mk I)$   
we take a Kuranishi chart $\mc K_{\gamma}=(\wh{\mb M}_{\gamma}, G_{\gamma}, \mu_{\gamma})$ 
naturally obtained 
by a closed subset $\wh{\mb M}_{\gamma}$ of $\wh{\mb M}$ and a closed subgroup 
$G_{\gamma}$ of $G$.
We show that we have a natural isomorphism
\begin{align}
\label{ambient}
[\wh{\mb M}_{\gamma}^{ss}/G_{\gamma}] \cong \mc N_{\gamma}
\end{align} 
where $\wh{\mb M}_{\gamma}^{ss}=\wh{\mb M}_{\gamma} \times_{\wh{\mb M}} \wh{\mb M}^{ss}$.
Furthermore we describe $\mu_{\gamma}$ in terms of moment map for the framed quivers appeared in
the preceding sections.

We also compute {\it virtual normal bundles} $N_{\gamma}$ of $\mc M_{\gamma}$ in 
$\mc M$ by $\mk N(\mc M_{\gamma})$.
By \eqref{localDM}, we have 
\begin{align}
\label{local}
[\mc M]^{vir} = \iota_{\ast} \sum_{\gamma=\pm, \mk I^{\sharp} \in \mc D^{\ell}(\mk I)} 
\frac{[\mc M_{\gamma}]^{vir}}{\Eu(N_{\gamma})}  
\end{align}
in $A^{\C^{\ast}_{\hbar}\times \mb{T}}_{\bullet}( \mc M ) \otimes \mb Q (\!( \hbar^{-1} )\!)$,
where $\iota \colon \mc M^{\C^{\ast}_{\hbar}} \to \mc M$ is the inclusion.
\begin{NB}
We have $[\mc N] = j_{\ast} \sum \frac{[\mc N_{\gamma}]}{e(\mk N(\mc N_{\gamma}))}$ for a smooth DM stack with $\C^{\ast}$-action, where $j \colon \mc N^{\C^{\ast}_{\hbar}} \to \mc N$ is the embedding.
For the diagram
\[
\xymatrix{
\mc N_{\gamma} \ar[r]^{j_{\gamma}} & \mc N \ar[r]^{e} & \E \\
\ar[u] \mc M_{\gamma}^{\iota_{\gamma}} \ar[r] &\ar[u] \mc M \ar[r] & \mc N \ar[u]^{0}},
\]
we have $0^{!} j_{\gamma \ast} = \iota_{\gamma \ast} 0^{!}$, where $0 \colon \mc N \to \E$ is the zero-section.
Furthermore for the diagram
\[
\xymatrix{
\mc N_{\gamma} \ar[r]^{j_{\gamma}} & \E_{\gamma} \ar[r] & E_{\gamma} \oplus (\E|_{\mc N_{\gamma}})^{mov} \\
\ar[u] \mc M_{\gamma}^{\iota_{\gamma}} \ar[r] &\ar[u]^{0_{\gamma}} \mc N_{\gamma} \ar@{=}[r] & \mc N_{\gamma} \ar[u]^{0|_{\mc N_{\gamma}}}},
\]
we have $0^{!}= e((\E|_{\mc N_{\gamma}})^{mov}) \cap 0_{\gamma}^{!}$ by the excess intersection formula, where $0_{\gamma} \colon \mc N_{\gamma} \to \E_{\gamma}$ is the zero-section.
Hence applying the Gysin homomorphism $0^{!}$ to the equation $[\mc N] = j_{\ast} \sum \frac{[\mc N_{\gamma}]}{e(\mk N(\mc N_{\gamma}))}$, we get
\begin{align*}
[\mc M]^{vir} = \iota_{\ast} \sum \frac{[\mc M_{\gamma}]^{vir}}{e(\mk N(\mc M_{\gamma}))}. 
\end{align*}.
\end{NB}
In the following computations, we cite corresponding results in \cite[\S 5]{M}, but they are much simpler in our case.

\subsection{Kuranishi chart of $\mc M_{\pm}$ and normal bundle $N_{\pm}$}
\label{subsec:obst1}
For $\gamma=\pm$, we take $\wh{\mb M}_{\pm}=\Proj \text{Sym} (\mc L_{\pm})$ and $G_{\pm}=G$.
We see that $\C^{\ast}_{\hbar}$-action in \eqref{act0} is trivial on $E|_{\wh{\mb M}}$, where 
$E=\wh{\mb M} \times \mb L(V)$. 
Hence $N_{\pm} = N_{\mc N_{\alpha}}|_{\wh{\mb M}_{\pm}}$, and 
the following lemma holds.
\begin{lem}[\protect{\cite[Proposition 5.9.2]{M}}]
\label{ob1}
We have 
$N_{\pm} = \mc L_{\mp} \otimes \mc L_{\pm}^{\vee} \otimes e^{\pm \hbar}$, 
where $\mc L_{\pm}$ are the induced line bundles on $\mc M_{\pm}$ from \eqref{linearization}.
\end{lem}
The tautological bundle $\mo (1)$ on $\wh{\mb M}_{\pm}=\Proj \text{Sym} (\mc L_{\pm})$
gives the $G$-polarization and $\mu \colon \wh{\mb M}_{\pm} 
= \wt{\mb M} \to \mb L_{\tilde{Q}}(V)$ is the moment map for $\tilde{Q}$.

\subsection{Kuranishi chart of $\mc M_{\mk I^{\sharp}}$}
\label{subsec:obst2}
For $\gamma=\mk I^{\sharp} \in \mc D^{\ell}(\mk I)$, we fix 
a direct sum decomposition $V= V^{\sharp} \oplus V^{\flat}$ 
and $\tilde{V}= \tilde{V}^{\sharp} \oplus \tilde{V}^{\flat}$ 
according to the decomposition type $\mk I^{\sharp}$
as in \S \ref{subsec:modified}.
This gives a closed embedding 
$\mb M_{\tilde{Q}}(\tilde{V}^{\sharp}) \times \mb M_{\tilde{Q}}( \tilde{V}^{\flat}) \to 
\wt{\mb M} = \mb M_{\tilde{Q}}(\tilde{V})$.
We put 
\[
\wh{\mb M}_{\mk I^{\sharp}} = \left(\mb M_{\tilde{Q}}(\tilde{V}^{\sharp}) \times \mb M_{\tilde{Q}}( \tilde{V}^{\flat}) \right) 
\times_{\wt{\mb M}} \wh{\mb M}
\] 
and $G_{\mk I^{\sharp}}= G_{\tilde{V}^{\sharp}} \times G_{\tilde{V}^{\flat}}$.
Then we have an isomorphism \eqref{ambient} since elements in 
$\mc N_{\mk I^{\sharp}}$ are included in the image of 
$\wh{\mb M}_{\mk I^{\sharp}}$ up to $G$-actions, and their stabilizers 
are included in $\C^{\ast} \id_{V^{\sharp}} \subset G_{\mk I^{\sharp}}$ 
by Lemma \ref{lem:stabilizer}.

Furthermore taking finite covers of $\C^{\ast}_{\hbar}$, 
we can modify $\C^{\ast}_{\hbar}$-actions up to $G_{\mk I^{\sharp}}$-actions 
as in \eqref{act} 
such that $\C^{\ast}_{\hbar}$ trivially acts on $\wt{\mb M}_{\mk I^{\sharp}}$. 
Hence we have decompositions $E|_{\wh{\mb M}_{\mk I^{\sharp}}}=
E_{\mk I^{\sharp}} \oplus (E|_{\wh{\mb M}_{\mk I^{\sharp}}})^{mov}$ by invariant parts 
$E|_{\mk I^{\sharp}}$ and moving parts $(E|_{\wh{\mb M}_{\mk I^{\sharp}}})^{mov}$ with respect to 
$\C^{\ast}_{\hbar}$-action where $E = \wh{\mb M} \times \mb L_{Q}(V)$.
Since the section $\sigma=\id_{\wh{\mb M}} \times \hat{\mu}$ is $\C^{\ast}_{\hbar}$-invariant, 
we see that $\sigma_{\mk I^{\sharp}}= \sigma|_{\wh{\mb M}_{\mk I^{\sharp}}}$ belongs to the invariant part 
$\Gamma(\wh{\mb M}_{\mk I^{\sharp}}, E_{\mk I^{\sharp}})$.

We take a Kuranishi chart $\mc K_{\mk I^{\sharp}}=
(\wh{\mb M}_{\mk I^{\sharp}}, G_{\mk I^{\sharp}}, \sigma_{\mk I^{\sharp}})$ of 
$\mc M_{\mk I^{\sharp}}$.
The semi-stable locus in $\wh{\mb M}_{\mk I^{\sharp}}$ is contained in
$\mb M_{\tilde{Q}}(\tilde{V}^{\sharp}) \times \mb M_{\tilde{Q}}(\tilde{V}^{\flat}) 
\times \C^{\ast} \chi_{\theta^{+}} / \chi_{\theta^{-}}$.
Hence we can use a description in \S \ref{subsec:modified}.
Furthermore we have $E_{\mk I^{\sharp}}= 
\wh{\mb M}_{\mk I^{\sharp}} \times \mb L_{Q}(V^{\sharp}) \times \mb L_{Q}(V^{\flat})$, and
$\sigma_{\mk I^{\sharp}}= \id_{\wh{\mb M}_{\mk I^{\sharp}}} \times \hat{\mu}_{\mk I^{\sharp}}$ 
where $\hat{\mu}_{\mk I^{\sharp}} \colon \wh{\mb M}_{\mk I^{\sharp}} \to 
\mb L_{Q}(V^{\sharp}) \times \mb L_{Q}(V^{\flat})$ is the composition of the projection to 
$\mb M_{\tilde{Q}}(\tilde{V}^{\sharp}) \times \mb M_{\tilde{Q}}(\tilde{V}^{\flat})$ and $\tilde{\mu} \times \tilde{\mu}$.
By Lemma \ref{lem:flatsharp} and Lemma \ref{lem:sharp}, 
the restriction $\mo(1)|_{\wh{\mb M}_{\mk I^{\sharp}}}$ gives 
$(\wh{\mb M}_{\mk I^{\sharp}})^{ss}=
\wt{\mb M}_{\tilde{Q}}^{\zeta^{\sharp}} (\tilde{V}^{\sharp})
\times \wt{\mb M}_{\tilde{Q}}^{\bar{\zeta}, \min(\mk I^{\sharp})-1} 
(\tilde{V}^{\flat}) \times \C^{\ast} \chi_{\theta^{+}} / \chi_{\theta^{-}}$. 

\subsection{Kuranishi chart of decompositions}
\label{subsec:obst3}
To get a Kuranishi chart $\mc K'_{\mk I^{\sharp}}=(\wh{\mb M}_{\mk I^{\sharp}}', G_{\mk I^{\sharp}}', 
\sigma_{\mk I^{\sharp}}')$ of 
$\mc M_{\mk I^{\sharp}}'$, we note that the description \eqref{decomp0} of an open embedding 
$\mc M_{\mk I^{\sharp}} \to \mc M$ is obtained by a natural coordinate change of $\mc K$.
We put 
\[
\wh{\mb M}_{\mk I^{\sharp}}'=\mb M_{\tilde{Q}}(\tilde{V}^{\sharp}) \times \C^{\ast} \chi_{\sharp} t^{-1/d^{\sharp}D}  
\times \mb M_{\tilde{Q}}(\tilde{V}^{\flat}) \times \C^{\ast} \chi_{\flat} t^{1/d^{\sharp} D}, 
\quad
G_{\mk I^{\sharp}}'=
G_{\tilde{V}^{\sharp}} \times G_{\tilde{V}^{\flat}} \times \C^{\ast}_{t^{1/d^{\sharp} D}}
\]
as in \S \ref{subsec:modified}.
We have a natural map $\wh{\mb M}_{\mk I^{\sharp}}' \to \wh{\mb M}_{\mk I^{\sharp}}$ induced by 
a multiplication $\chi^{\sharp} \chi^{\flat}  =\left. (\chi_{\theta^{+}} / \chi_{\theta^{-}} ) \right|_{G_{\mk I^{\sharp}}}$, 
where $G_{\mk I^{\sharp}} = G_{\tilde{V}^{\sharp}} \times G_{\tilde{V}^{\flat}}$ is regarded as a sub-group
of $G_{\tilde{V}}$ via the decomposition $\tilde{V}=\tilde{V}^{\sharp} \oplus \tilde{V}^{\flat}$.
$G_{\mk I^{\sharp}}'$-invariant section $\sigma_{\mk I^{\sharp}}'$ and $G_{\mk I^{\sharp}}'$-linearization 
are obtained via pull-back by this map.
The \'etale morphism $\Phi_{\mk I^{\sharp}}$ is obtained from the coordinate change of 
Kuranishi charts $\mc K_{\mk I^{\sharp}}$ and $\mc K'_{\mk I^{\sharp}}$ 
induced by these maps $\wh{\mb M}_{\mk I^{\sharp}}' \to \wh{\mb M}_{\mk I^{\sharp}}$ 
and $G_{\mk I^{\sharp}}' \to G_{\mk I^{\sharp}}$. 

We also take Kuranishi charts 
$\mc K_{\sharp}=(\mb M_{\sharp},  G_{\sharp}, \varphi_{\sharp})$ and 
$\mc K_{\flat}=(\mb M_{\flat}, G_{\flat}, \varphi_{\flat})$ 
of $Fl_{H_{Q}(d^{\sharp} \beta)}( \bar{\mc V}^{\sharp}_{0}, \bar{\mk I}^{\sharp})$ and 
$\wt M^{\min(\mk I^{\sharp})-1}(\alpha - d^{\sharp} \beta, \mk I^{\flat})$ as follows.
For $\mc K_{\sharp}$, we put $\mb M_{\sharp}= \mb M_{\tilde{Q}}(\tilde{V}^{\sharp})$,  
\[
G_{\sharp} = \prod_{i \in I} \GL(V^{\sharp}_{i}) \times \prod_{k > \min(\mk I^{\sharp})} \GL(V^{\sharp}_{(0,k)}),
\] and 
$\varphi_{\sharp} = \tilde{\mu}$.
For $\mc K_{\flat}$, we put
$\mb M_{ \flat} =\mb M_{\tilde{Q}}(\tilde{V}^{\flat}), G_{\flat}=G_{\tilde{V}^{\flat}}$, and 
$\varphi_{\flat} =\tilde{\mu}$ as usual notation.
By the argument in the proof of Theorem \ref{decomp},
the \'etale covering
\[
\Psi_{\mk I^{\sharp}} \colon \mc M_{\mk I^{\sharp}}' \to  
Fl_{H_{Q}(d^{\sharp} \beta)}( \bar{\mc V}^{\sharp}_{0}, \bar{\mk I}^{\sharp}) \times 
\wt M^{\min(\mk I^{\sharp})-1}(\alpha - d^{\sharp} \beta, \mk I^{\flat})
\]
is obtained by coordinate change of $\mc K_{\mk I^{\sharp}}'$ and 
$\mc K_{\sharp} \times \mc K_{\flat} =(\mb M_{\sharp} \times \mb M_{\flat},  
G_{\sharp} \times G_{\flat}, \varphi_{\sharp} \times \varphi_{\flat})$.

\subsection{Normal bundle $N_{\mk I^{\sharp}}$}
\label{subsec:norm}
In the following, for vector bundles $\mc E, \mc F$ on a stack $\mc Z$, we write by 
$\mc Hom(\mc E, \mc F)$ the vector bundle $\mc E^{\vee} \otimes \mc F$ on $\mc Z$.

We consider $\C^{\ast}_{e^{\hbar/d^{\sharp}D}}$-action on $\mc M$ as in \eqref{act}. 
We take universal full flags $\bar{\mc F}_{\bullet}^{\sharp}$ on 
$Fl_{H_{Q} (d^{\sharp} \beta )}( \bar{\mc V}^{\sharp}_{0}, \bar{\mk I}^{\sharp})$, 
and $\mc F_{\bullet}^{\flat}$ on $M_{\tilde{Q}}^{\bar{\zeta}, \min(\mk I_{\sharp}) -1} (\tilde{V}_{\flat})$.
On $Fl_{H_{Q} (d^{\sharp} \beta )}( \bar{\mc V}^{\sharp}_{0}, \bar{\mk I}^{\sharp})$,
we define a full flag $\mc F_{\bullet}^{\sharp}$ of $\mc V^{\sharp}_{0}$ 
by the pull-back of $\bar{\mc F}_{\bullet}^{\sharp}$ to $\mc V^{\sharp}_{0}$ for $i \neq \min(\mk I^{\sharp})$
and $\mc F_{\min(\mk I^{\sharp})}^{\sharp} = \im (\mc V^{\sharp}_{\infty'} \to \mc V^{\sharp}_{0})$. 
We also write their pull-backs to the product $\mc M_{\sharp} \times \mc M_{\flat}$
by the same letter $\mc F^{\sharp}_{\bullet}, \mc F^{\flat}_{\bullet}$.

For two flags $\mc F_{\bullet}, \mc F_{\bullet}'$ of sheaves on the same Deligne-Mumford stack,  
we put 
\begin{align*}
\Theta (\mc F_{\bullet}, \mc F_{\bullet}') 
&=
\sum_{i > j }  
\mc Hom\left(\mc F_{j} / \mc F_{j-1}, \mc F_{i}' /\mc F_{i-1}' \right),
\\
\mk H (\mc F_{\bullet}, \mc F_{\bullet}')
&=
\Theta (\mc F_{\bullet}, \mc F_{\bullet}') 
+
\Theta (\mc F_{\bullet}', \mc F_{\bullet}).
\end{align*}
When $\mc F_{\bullet} = \mc F_{\bullet}'$, we put $\Theta(\mc F_{\bullet}) = \Theta(\mc F_{\bullet}, \mc F_{\bullet})$.

\begin{lem} 
\label{nmki}
We have $\Phi_{\mk I^{\sharp}}^{\ast} N_{\mk I^{\sharp} } \cong 
\Psi_{\mk I^{\sharp}}^{\ast} \left[ \mathfrak{N}( \mc V^{\sharp} \otimes
e^{\hbar/d^{\sharp}D} \otimes 
(L_{\mk I^{\sharp}})^{\vee}, \mc V^{\flat})  + 
\mk H (\mc F^{\sharp}_{\bullet}, \mc F^{\flat}_{\bullet}) \right]
$ where
\begin{align*}
\mathfrak{N}( \mc V^{\sharp}, \mc V^{\flat})
&=
\sum_{(\heartsuit, \spadesuit)} 
\left(
\sum_{\substack{a \in Q_{1} }} 
\mc H om( \mc V^{\heartsuit}_{\out(a)}, \mc V^{\spadesuit}_{\inn(a)}) 
- 
\sum_{l \in Q_{2}} \mc H om( \mc V^{\heartsuit}_{\out(l)}, \mc V^{\spadesuit}_{\inn(\l)}) 
-
\sum_{i \in I} 
\mc Hom ( \mc V^{\heartsuit}_{i}, \mc V^{\spadesuit}_{i}) \right),
\notag
\end{align*}
where the sum is taken for $(\heartsuit, \spadesuit)= (\sharp, \flat), (\flat, \sharp)$.
\end{lem}
\proof
In the $K$-group of $\mc M_{\mk I^{\sharp}}'$, we have
\begin{align*}
&
\sum_{(\heartsuit, \spadesuit)} 
\left(
\sum_{a \in \tilde{Q}_{1} \setminus Q_{1}}  
\mc H om( \mc V^{\heartsuit}_{\out(a)}, \mc V^{\spadesuit}_{\inn(a)}) -
\sum_{v \in \tilde{Q}_{0} \setminus Q_{0}} \mc Hom ( \mc V^{\heartsuit}_{v}, \mc V^{\spadesuit}_{v}) 
\right) \\
& =
\sum_{(\heartsuit, \spadesuit)} \sum_{i > j }  
\mc Hom
\left(\mc F_{j}^{\heartsuit} / \mc F_{j-1}^{\heartsuit}, \mc F_{i}^{\spadesuit} /\mc F_{i-1}^{\spadesuit} \right)
=\mk H (\mc F_{\bullet}^{\heartsuit}, \mc F_{\bullet}^{\spadesuit}),
\end{align*} 
where the sum is also taken for $(\heartsuit, \spadesuit)=(\flat, \sharp), (\sharp, \flat)$.
\endproof 

\subsection{Relative tangent bundles for flags}
 \label{subsec:rela2}
We consider the pull-back $\Theta^{rel}$ to $\mc M$ of the relative tangent bundle of $[\wt{\mu}^{-1}(0)/G]$ over $[\mu^{-1}(0)/G]$. 
After restricting to $\mc M_{\mk I^{\sharp}}$ and pulling back to $\mc M_{\mk I^{\sharp}}'$, we have a decomposition 
\begin{eqnarray}
\label{rel}
\Phi_{\mk I^{\sharp}}^{\ast} \Theta^{rel} = \Psi_{\mk I^{\sharp}}^{\ast} 
( \Theta(\mc F_{\bullet}^{\sharp}) \oplus 
\Theta(\mc F_{\bullet}^{\flat}) \oplus \mk H (\mc F_{\bullet}^{\sharp}, \mc F_{\bullet}^{\flat})  )
\end{eqnarray} 
We also consider the relative tangent bundle $\Theta'$ of 
the projection 
$Fl_{H_{Q} (d^{\sharp} \beta)}( \bar{\mc V}^{\sharp}_{0}, \bar{\mk I}^{\sharp}) \to
H_{Q} (d^{\sharp} \beta)$.
Then we have an exact sequence
\begin{align}
\label{ob4}
0 \to \Theta' \to \Theta(\mc F_{\bullet}^{\sharp}) \to \bar{\mc V}^{\sharp}_{0} \to 0,
\end{align}
where $\bar{\mc V}^{\sharp}_{0}$ is the quotient 
by the tautological homomorphism $\mc V^{\sharp}_{\infty'} \to \mc V^{\sharp}_{0}$, 
and it is identified with $\mc Hom (\mc F^{\sharp}_{\min(\mk I^{\sharp})}, 
\mc V^{\sharp}_{0}/ \mc F^{\sharp}_{\min(\mk I^{\sharp})})$.

\begin{NB}
For any full flag bundle $p \colon Y =Fl(\E, \underbar{n}) \to X$ 
for a rank $n$ vector bundle $\E$ on $X$, 
we can check $p_{\ast} \left( e(T_{Y/X}) \cap p^{\ast} \beta \right) = 
n! \beta$ for $\beta \in A_{\bullet}(X)$ as follows.
Since full flag bundles are obtained by iterations of projective bundles, 
it is enough to show that for $\pi \colon Z =\PP(\E, \underbar{n}) \to X$, 
we have $\pi_{\ast} \left( e(T_{Z/X}) \cap \pi^{\ast} \beta \right)  = n \beta$ for $\beta \in A_{\bullet}(X)$.
This follows from the description of $\pi_{\ast} \colon A_{\bullet}(Z) \to A_{\bullet}(X)$ in \cite{F}.
\end{NB}

\subsection{Localization}

We put $M_{0} = \Spec \Gamma(\mu^{-1}(0), \mo_{\mu^{-1}(0)})^{\prod_{i \in I} \GL(V_{i})}$.
We consider the naturally induced proper morphisms $\Pi \colon \mc M \to M_{0}$.
We also write by the same letter $\Pi$ the proper morphisms from various moduli stack $M$ to $M_{0}$. 
We assume that $\Pi \colon M \to M_{0}$ is $\mb T \times \C^{\ast}_{\hbar}$-equivariant
and have Kuranishi charts compatible with $\mb T \times \C^{\ast}_{\hbar}$-action
as in \S \ref{subsec:virt} and \S \ref{subsec:local1}.
Here $\C^{\ast}_{\hbar}$ acts on $M_{0}$ trivially.

For $\varphi$ in $\mb T \times \C^{\ast}_{\hbar}$-equivariant Chow ring 
$A^{\bullet}_{\mb T \times \C^{\ast}_{\hbar}}$, we write by $\int_{M} \varphi$ the push-forward
$\Pi_{\ast} (\varphi \cap [M]^{vir})$ in the Chow group 
$A^{\mb T \times \C^{\ast}_{\hbar}}_{\bullet}(M_{0}) \otimes_{A} S$.
Here $A=A^{\bullet}_{\mb T \times \C^{\ast}_{\hbar}}(\text{pt})$ is 
the $\mb T \times \C^{\ast}_{\hbar}$-equivariant Chow ring
of the point, and $S$ is the fractional field of $A$.  
In the following, we assume that the fixed points set $(M_{0})^{\mb T}$ consists of one point. 
We identify $A^{\mb T \times \C^{\ast}_{\hbar}}_{\bullet}(M_{0}) \otimes_{A} S$
and $A^{\mb T \times \C^{\ast}_{\hbar}}_{\bullet}(M_{0}^{\mb T}) \otimes_{A} S \cong S$ 
via the isomorphism $\iota_{\ast} \otimes_{A} S$ by the embedding 
$\iota \colon (M_{0} )^{\mb T} \to M_{0}$.

By Theorem \ref{decomp} and \eqref{local}, we have the following diagram :
\[
\xymatrix{
A^{\bullet}_{\C^{\ast}_{\hbar}\times \mb{T}} (\mc M) \otimes_{\C[\hbar]} \C(\!(\hbar^{-1})\!) 
\ar[d]_{\int_{[\mc M]^{vir}} } \ar[r]^-{\cong}& 
A^{\bullet}_{\mb{T}} (\mc M^{\C^{\ast}_{\hbar}}) \otimes_{\C[\hbar]} \C(\!(\hbar^{-1})\!)
\ar[d]^{\int_{ [\mc M_{+}]^{vir}} +  
\int_{[\mc M_{-}]^{vir}} + \sum_{\mk I^{\sharp}} \int_{[\mc M_{\mk I^{\sharp}}]^{vir}} } \\
A_{\bullet}^{\mb T} (M_{0}) \otimes_{\C[\hbar]} \C(\!(\hbar^{-1})\!) \ar@{=}[r] & 
A_{\bullet}^{\mb T} (M_{0}) \otimes_{\C[\hbar]} \C(\!(\hbar^{-1})\!)\\
}
\]
Here $\Eu(N_{\pm}), \Eu(N_{\mk I^{\sharp}})$ are invertible in 
$A^{\bullet}_{\tilde{\mb T} \times \C^{\ast}_{\hbar}}(\mc M)( \!( \hbar^{-1} )\!)$, and
for $\varphi \in A^{\bullet}_{\tilde{\mb T} \times \C^{\ast}_{\hbar}}(\mc M)$.
The upper horizontal arrow is given by 
$\frac{\iota_{+}^{\ast}}{\Eu (N_{+})} + \frac{\iota_{-}^{\ast}}{\Eu (N_{-})} + 
\sum_{\mk I^{\sharp} \in \mc D^{\ell}(\mk I)} \frac{\iota_{\mk I^{\sharp}}^{\ast}}{\Eu (N_{\mk I^{\sharp}})}$.

We write by $\hbar$ the first Chern class of the weight space $\C_{e^{\hbar}}$ in 
$A^{\bullet}_{\C^{\ast}_{\hbar}}(\text{pt})$ of the weight $e^{\hbar} \in \C^{\ast}_{\hbar}$, and $\iota_{\pm}$ and $\iota_{\mk I^{\sharp}}$ are embeddings of $\mc M_{\pm}$ and $\mc M_{\mk I^{\sharp}}$ into $\mc M$.
We have
\begin{align}
\label{localization}
\int_{\mc M} \varphi = \int_{\mc M_{+}} \frac{\varphi|_{\mc M_{+}}}{\Eu (N_{+})} + \int_{\mc M_{-}} \frac{\varphi|_{\mc M_{-}}}{\Eu (N_{-})} + \sum_{\mk I^{\sharp}} \int_{\mc M_{\mk I^{\sharp}}} \frac{\varphi|_{\mc M_{\mk I^{\sharp}}}}{\Eu (N_{\mk I^{\sharp}})}.
\end{align}

We introduce another torus $\C^{\ast}_{e^{\theta}}$ trivially acting on all moduli spaces.
For $K$-theory class $\alpha=E-F$, we put $\Eu^{\theta} (\alpha) = \Eu(E \otimes \C_{e^{\theta}})/
\Eu(F \otimes \C_{e^{\theta}})$.
We take the cohomology class $\psi \in A_{\C^{\ast}_{\hbar} \times \mb T}^{\bullet}( \mc M)$.
Furthermore, we put 
\begin{eqnarray*}
\tilde{\psi} =  \frac{ \psi \cdot \Eu^{\theta} (\Theta^{rel})}{|\mk I| !} \in 
A_{\C^{\ast}_{\hbar} \times \mb T}^{\bullet}( \mc M )  
\end{eqnarray*}
and substitute $\varphi = \tilde{\psi}$ in \eqref{localization}.
The left hand side in \eqref{localization} is a polynomial in $\hbar$ while the right hand side has a power series part in $\hbar^{-1}$.
Hence if the symbol $\displaystyle \Res_{\hbar = \infty}$ denotes the operation taking the coefficient in $\hbar^{-1}$, we have
\begin{align*}
\int_{\wt M^{\ell}(\alpha, \mk I)} \tilde{\psi}
-\int_{M^{\zeta^{-}}(\alpha)} \psi 
\notag
&=
- \Res_{\hbar=\infty} \sum_{\mk I^{\sharp} \in D^{\ell}(\mk I)} \int_{\mc M_{\mk I^{\sharp}}} 
\frac{\tilde{\psi}|_{\mc M_{\mk I^{\sharp}}}}{\Eu (N_{\mk I^{\sharp}})}.
\end{align*}

By Lemma \ref{nmki} and \eqref{rel},
the last summand is equal to
\begin{align}
\label{euler}
\frac{|\mk I^{\flat}|!}{|\mk I|!} 
\int_{\wt M^{ \min(\mk I^{\sharp})-1} (\alpha-  d^{\sharp} \beta, \mk I^{\flat})} 
\frac{1}{|\mk I^{\flat}|!} 
\int_{Fl_{H_{Q}(d^{\sharp} \beta)}( \bar{\mc V}^{\sharp}_{0}, \bar{\mk I}^{\sharp})} 
\frac{
\psi \cdot
\Eu^{\theta}(\Theta( \mc F^{\sharp}_{\bullet} \oplus \mc F^{\flat}_{\bullet} ))
}
{\Eu(\mathfrak{N}(\mc V^{\sharp} \otimes e^{\hbar}, \mc V^{\flat})) \cdot 
\Eu(\mk H (\mc F^{\sharp}_{\bullet}, \mc F^{\flat}_{\bullet}))
}.
\end{align}
Here we have a proper morphism from
\[
\Pi^{\sharp} \colon
\wt M^{ \min(\mk I^{\sharp})-1} (\alpha-  d^{\sharp} \beta, \mk I^{\flat}) \times
Fl_{H_{Q}(d^{\sharp} \beta)}( \bar{\mc V}^{\sharp}_{0}, \bar{\mk I}^{\sharp}) \to
\wt M^{ \min(\mk I^{\sharp})-1} (\alpha-  d^{\sharp} \beta, \mk I^{\flat})
 \times M_{0}(V^{\sharp} \oplus V_{\infty'}).
 \]
In this situation, we write by $\int_{Fl_{H_{Q}(d^{\sharp} \beta)}( \bar{\mc V}^{\sharp}_{0}, \bar{\mk I}^{\sharp})}$
the push-forward $\left( \Pi^{\sharp} \right)_{\ast} \varphi \cap
[Fl_{H_{Q}(d^{\sharp} \beta)}( \bar{\mc V}^{\sharp}_{0}, \bar{\mk I}^{\sharp})]$
in $A_{\bullet}^{\mb T \times \C^{\ast}_{\hbar}}
(\wt M^{ \min(\mk I^{\sharp})-1} (\alpha-  d^{\sharp} \beta, \mk I^{\flat}) 
\times M_{0}(V^{\sharp} \oplus V_{\infty'}))
\otimes_{A} S \cong A_{\bullet}^{\mb T \times \C^{\ast}_{\hbar}}
(\wt M^{ \min(\mk I^{\sharp})-1} (\alpha-  d^{\sharp} \beta, \mk I^{\flat}) )
\otimes_{A} S$
for $\varphi \in A^{\bullet}_{\mb T}
(Fl_{H_{Q}(d^{\sharp} \beta)}( \bar{\mc V}^{\sharp}_{0}, \bar{\mk I}^{\sharp}))$. 

In the expression \eqref{euler}, we deleted some line bundles and a parameter $d^{\sharp} D$, since we have 
$\displaystyle \Res_{\hbar = \infty} f(\hbar) = 
d^{\sharp} D \Res_{\hbar = \infty} f( d^{\sharp} D \hbar+a)$ (cf. \cite[\S 8.2]{O1}), 
and integrals over degree $1/(d^{\sharp} D)$ \'etale covering of full flag bundles
$Fl_{H_{Q}(d^{\sharp} \beta)}( \bar{\mc V}^{\sharp}_{0}, \bar{\mk I}^{\sharp})$.

\section{Wall-crossing formula}
In the following, we deduce wall-crossing formula from analysis in the previous section.
These are the similar calculations to \cite[\S 6]{NY2}, hence we omit detail explanation.

In the following, we fix a $Q_{0}$-graded vector space $V=\bigoplus_{v \in Q_{0}} V_{v}$ 
with $\dim V_{\infty}=1$.
We put $\alpha=(\alpha_{v})_{v\in Q_{0}} = (\dim V_{v})_{v \in Q_{0}} \in (\Z_{\ge 0})^{Q_{0}}$.
For $\mbi d= (d_{1}, \ldots, d_{j}) \in \Z_{>0} ^{j}$, we put $|\mbi d| = d_{1} + \cdots + d_{j}$.
Let $\Dec{}_{\beta_{0}, j}^{\alpha_{0}}$ be the set of collections 
$\mbi {\mk I}=( \mk I_{1}, \ldots, \mk I_{j})$ of disjoint non-empty subsets of $[\alpha_{0}]=\lbrace 1, 2, \ldots, \alpha_{0} \rbrace$ 
such that
\begin{enumerate}
\item[$\bullet$] $|\mk I_{i}| = d_{i} \beta_{0}$ for $d_{i} \in \Z_{>0}$
and $i=1, \ldots,j$, and
\item[$\bullet$] $\min(\mk I_{1}) > \cdots > \min(\mk I_{j})$.
\end{enumerate}
We put $\mbi d_{\mbi{\mk I}} = 
\frac{1}{\beta_{0}}(|\mk I_{1}|, \ldots, |\mk I_{j}|)
=(d_{1}, \ldots, d_{j}) \in \Z_{>0}^{j}$
and $\mk I_{\infty} =[\alpha_{0}] \setminus \bigsqcup_{i=1}^{j} \mk I_{i}$.
We note that $\Dec{}_{\beta_{0}, 1}^{\alpha_{0}} 
= \mc D^{\ell}([\alpha_{0}])$
and
\[
\Dec{}_{\beta_{0}, j+1}^{\alpha_{0}}
=
\lbrace
(\mk I_{1}, \ldots, \mk I_{j}, \mk I_{j+1})
\mid 
\mbi{\mk I}=(\mk I_{1}, \ldots, \mk I_{j}) \in 
\Dec{}_{\beta_{0}, j}^{\alpha_{0}}, \ 
\mk I_{j+1} \in \mc D^{\min(\mk I_{j})-1}( \mk I_{\infty})
\rbrace.
\]
We define $\sigma \colon \Dec{}_{\beta_{0}, j+1}^{\alpha_{0}} \to \Dec{}_{\beta_{0}, j}^{\alpha_{0}}$ by
$\sigma(\mk I_{1}, \ldots, \mk I_{j}, \mk I_{j+1}) 
= (\mk I_{1}, \ldots, \mk I_{j})$.



\subsection{Iterated cohomology classes}
\label{subsec:iter1}
\indent For $\mbi d =(d_{1}, \ldots, d_{j})\in \Z_{>0}^{j}$ and $\ell=0, 1, \ldots, \alpha_{0}$, we consider a product 
\[
M_{\mbi d} = M^{\zeta^{-}} (\alpha - |\mbi d| \beta) 
\times \prod_{i=1}^{j} H_{Q}(d_{i} \beta)
\] 
We write by $\mc V^{(i)}$ the tautological bundle
$\mc V^{\sharp} \oplus \mc V^{\sharp}_{\infty'}$ on the component 
$H_{Q}(d_{i} \beta)$.

For $\mbi{ \mk I}=(\mk I_{1}, \ldots, \mk I_{j}) \in \rho^{-1}(\mbi d)$, we also put
\[
\wt M_{\mbi{\mk I}}^{\ell}  = \wt M^{\ell}(\alpha - |\mbi d_{\mbi{\mk I}}| \beta, \mk I_{\infty}) \times
Fl_{H_{Q}(d_{1} \beta)}(\bar{\mc V}_{0}^{(1)}, \bar{\mk I}_{1}) 
\times \cdots \times 
Fl_{H_{Q}(d_{j} \beta)}(\bar{\mc V}_{0}^{(j)}, \bar{\mk I}_{j})
\]
where $\bar{\mc V}_{0}^{(k)}=\mc V_{0}^{(k)} / \mc V^{(k)}_{\infty'}$, 
$\bar{\mk I}_{i}=\mk I_{i} \setminus \lbrace \min(\mk I_{i})\rbrace$,
and write by $\mc F_{\bullet}^{(i)}$ the pull-back to $\wt M_{\mbi{\mk I}}^{\ell}$ of 
the flag $\mc F_{\bullet}^{\sharp}$ on each component 
$Fl(\bar{\mc V}_{0}^{(i)}, \bar{\mk I}_{i})$.
We regard $\lbrace \mc V_{(0, k)} \rbrace_{k=1}^{\alpha_{0}}$
on $\wt M^{\ell}( \alpha - |\mbi d_{\mbi{\mk I}}| \beta, \mk I_{\infty} )$
as a flag $\mc F_{\bullet}$ of $\mc V_{0}$, and  
write by $\mc F_{\bullet}^{\infty}$ the pull-back to 
$\wt M_{\mbi{\mk I}}^{\ell}$.
We put $\mc F_{\bullet}^{>i}=\mc F_{\bullet}^{\infty} \oplus \bigoplus_{k > i} \mc F_{\bullet}^{(k)}$.
We also put $\mbi{\mk I}_{> i} = \mk I_{\infty} \sqcup \bigsqcup_{k >i} \mk I_{k}$.

For the fundamental cycle $[\wt M_{\mbi{\mk I}}^{\ell} ]$ defined from the obstruction theory in the previous section, 
and a cohomology class $\varphi \in A_{\mb {T} \times \C^{\ast}_{\theta} 
\times \prod_{i=1}^{j} \C^{\ast}_{\hbar_{i}}}^{\bullet}(\wt M_{\mbi{\mk I}}^{\ell} )$,  
we write by 
\[
\int_{[\wt M_{\mbi{\mk I}}^{\ell} ]} \varphi \in A_{\mb {T} \times \C^{\ast}_{\theta} \times 
\prod_{i=1}^{j} \C^{\ast}_{\hbar_{i}}}^{\bullet}
(\wt M^{\ell}(\alpha - |\mbi d_{\mbi{\mk I}}| \beta, \mk I_{\infty} ))
\]
the Poincare dual of the push-forward of $\varphi \cap [\wt M_{\mbi{\mk I}}^{\ell} ]$ 
by the projection $\wt M_{\mbi{\mk I}}^{\ell} \to \wt M^{\ell}(\alpha - |\mbi d_{\mbi{\mk I}}| \beta, \mk I_{\infty} )$. 

For the tautological bundle $\mc V$ on $\wt M^{\ell}(\alpha - |\mbi d_{\mbi{\mk I}}| \beta, \mk I_{\infty} )$, 
we write by the same letter 
the pull-backs to the product $\wt M_{\mbi{\mk I}}^{\ell}$.
For $\mc V^{(i)}$, we also write by the same letter the pull-back by the projection 
$\wt M_{\mbi{\mk I}}^{\ell} \to H_{Q}(d_{i} \beta)$
We put $\mc V^{>i} = \mc V \oplus \bigoplus_{k > i } \mc V^{(k)} \otimes e^{\hbar_{k}}$. 

For $\mbi{\mk I}=(\mk I_{1}, \ldots, \mk I_{j}) \in \text{Dec}_{\beta_{0}, j}^{\alpha_{0}} $, we write by 
$\tilde{\psi}_{\mbi{\mk I}} (\mc V)$
the following cohomology class
\begin{align}
&
\int_{[\wt M_{\mbi{\mk I}}^{\ell} ]} 
\frac{ 
\psi \left( \mc V \oplus \bigoplus_{i=1}^{j} \mc V^{(i)} \otimes e^{\hbar_{i}} \right) 
\Eu^{\theta}(\Theta(\mc F^{>0}_{\bullet}))
}
{ 
\Eu \left( \bigoplus_{i=1}^{j} \mathfrak{N} ( 
\mc V^{(i)} \otimes e^{\hbar_{i}}, 
\mc V^{>i} )
\right)
\prod_{i=1}^{j}
\Eu(\mk H(\mc F^{(i)}_{\bullet}, \mc F^{>i}_{\bullet}))
}
\frac{1}{|\mk I_{\infty}|!}
\label{itcoho}
\\
&=
\int_{[\wt M_{\mbi{\mk I}}^{\ell} ]} 
\frac{ 
\psi \left( \mc V \oplus \bigoplus_{i=1}^{j} \mc V^{(i)} \otimes e^{\hbar_{i}} \right) 
\Eu^{\theta}(\Theta(\mc F^{\infty}_{\bullet}))
/|\mk I_{\infty}|!
}
{ 
\Eu \left( \bigoplus_{i=1}^{j} \mathfrak{N} ( 
\mc V^{(i)} \otimes e^{\hbar_{i}}, 
\mc V^{>i} )
\right)
}
\prod_{i=1}^{j}
\Eu^{\theta}(\Theta(\mc F^{(i)}_{\bullet}))
\frac{
\Eu^{\theta} (\mk H ( \mc F^{(i)}_{\bullet}, \mc F^{>i}_{\bullet})) 
}
{
\Eu(\mk H(\mc F^{(i)}_{\bullet}, \mc F^{>i}_{\bullet}))
}
\notag
\end{align}
in $A^{\bullet}_{\mb{T} \times \C^{\ast}_{\theta} \times \prod_{i=1}^{j} \C^{\ast}_{\hbar_{i}}}(\wt M^{\ell}(\alpha - |\mbi d_{\mbi{\mk I}}| \beta ))$.
Here $e^{\hbar_{i}}$ is a trivial bundle with $e^{\hbar_{i}}$-weight.
By modified $\C^{\ast}_{\hbar}$-action \eqref{act}, we need to multiply $\mc V^{(i)}$ with $e^{\hbar_{i}}$ in \eqref{itcoho}.


\subsection{Recursions}
\label{subsec:loca}
For $\mbi{\mk I} = (\mk I_{1}, \ldots, \mk I_{j}) \in \Z_{>0}^{j}$, we put $\mk I = \mk I_{\infty}$,
$d=|\mbi d_{\mbi{\mk I}}|$, and $\ell=\min(\mk I_{j})-1$, and take 
an equivariant cohomology class $\varphi = \tilde{\psi}_{\mbi{\mk I }}(\mc V)$ on $\mc M$.
For the convenience, we also put $\tilde{\psi}_{()}=\psi$ for $j=0$. 
Then \eqref{euler} is equal to
\begin{eqnarray}
\label{exc}
\frac{|\mk I^{\flat}|! }{|\mk I| !} \int_{\wt M^{\min(\mk I^{\sharp})-1}(\alpha - d\beta - d^{\sharp} \beta ) } 
\tilde{\psi}_{ ( \mbi{\mk I}, \mk I^{\sharp} )}(\mc V),
\end{eqnarray}
where $d^{\sharp} = |\mk I^{\sharp}|$, and $(\mbi{\mk I}, \mk I^{\sharp})
=(\mk I_{1}, \ldots, \mk I_{j}, \mk I^{\sharp}) \in \Dec^{\alpha_{0}}_{\beta_{0}, j+1}$.

Using $\tilde{\psi}_{\mbi{\mk I}} (\mc V)$ defined in \eqref{itcoho},
we deduce recursion formula.
\begin{lem}
For $l \ge 1$, we have
\begin{align}
\notag
&
\int_{M^{\zeta^{+}}(\alpha)} \psi(\mc V)  - \int_{M^{\zeta^{-}}(\alpha )} \psi(\mc V) \\
\notag
&= 
\sum_{j=1}^{l-1} 
\Res_{\hbar_{1} = \infty} \cdots \Res_{\hbar_{j}=\infty}
\sum_{\mbi{\mk I} \in \Dec{}_{\beta_{0}, j}^{\alpha_{0}}} \frac{ |\mk I_{\infty}|! }{ \alpha_{0}! }
\int_{\wt M^{ 0}(\alpha - |\mbi d_{\mbi{\mk I}}|\beta)} \tilde{\psi}_{\mbi{\mk I}} (\mc V)\\
\label{formula2>0}
&+
\Res_{\hbar_{1} = \infty} \cdots \Res_{\hbar_{l}=\infty}
\sum_{\mbi{\mk I} \in \Dec{}_{\beta_{0}, l}^{\alpha_{0}}} \frac{|\mk I_{\infty}|! }{ \alpha_{0}! }
\int_{\wt M^{\min(\mk I_{l})-1} (\alpha - |\mbi d_{\mbi{\mk I}}| \beta )} \tilde{\psi}_{\mbi{\mk I}} (\mc V).
\end{align}
\end{lem}
\proof
We prove by induction on $j$.
For $l=1$, \eqref{formula2>0} is nothing but \eqref{exc} for $j_{0}=0$ and $\ell = n$.
For $l \ge 1$, we assume the formulas \eqref{formula2>0}.
Then again by \eqref{exc}, the last summand for each $\mbi{\mk I} \in \Dec{}_{\beta_{0}, l}^{\alpha_{0}}$ is equal to 
\begin{align*}
&
 \frac{|\mk I_{\infty}|!  }{ \alpha_{0}! }
\left( \int_{\wt M^{0} (\alpha - |\mbi d_{\mbi{\mk I}}| \beta )}  \tilde{\psi}_{\mbi{\mk I}} (\mc V)  \right. \\
&+
\Res_{\hbar_{l+1}=\infty}
\left. \sum_{\mbi{\mk I}' \in \sigma^{-1} (\mbi{\mk I})} 
\frac{|{\mbi{\mk I}'}^{c}|! }{|\mk I_{\infty}|!} 
\int_{\wt M^{\min (\mk I_{l+1}' )-1 }(\alpha - |\mbi d_{\mbi{\mk I}'}| \beta )} \tilde{\psi}_{\mbi{\mk I}'} (\mc V)  \right ),
\end{align*}
where $\mbi{\mk I}'=(\mbi{\mk I}, \mk I_{l+1}')$.
Hence we have \eqref{formula2>0} for general $l \ge 1$.
\endproof

For $l > \alpha_{0} / \beta_{0}$, the set $\Dec{}_{\beta_{0}, l}^{\alpha_{0}}$ is empty.
Thus we get the following theorem.
\begin{thm}
\label{thm:wcm}
We have
\begin{align}
&
\int_{M^{\zeta^{+}}(\alpha)} \psi(\mc V)  - \int_{M^{\zeta^{-}}(\alpha )} \psi(\mc V) 
\notag
\\
&=
\label{main>0}
\sum_{j=1}^{\lfloor \alpha_{0}/\beta_{0} \rfloor } \sum_{\mbi{\mk I} \in \Dec{}_{\beta_{0}, j}^{\alpha_{0}}} 
\frac{ |\mk I_{\infty}|!  }{ \alpha_{0}! }
\Res_{\hbar_{1} = \infty} \cdots \Res_{\hbar_{j}=\infty}
\int_{\wt M^{0} (\alpha - |\mbi d_{\mbi{\mk I}}| \beta )} \tilde{\psi}_{\mbi{\mk I}} (\mc V)
\end{align}
where $\tilde{\psi}_{\mbi{\mk I}} (\mc V)$ is defined in \eqref{itcoho}.
\end{thm}

\section{Euler class of the tangent bundle}
\label{sec:euler2}
For a quiver $Q=(Q_{0}, Q_{1}, Q_{2})$ and $Q_{0}$-graded vector bundle $\mc V$, we put 
\[
\Lambda_{Q}(\mc V) =
\sum_{a \in Q_{1}} \mc Hom (\mc V_{\out (a)}, \mc V_{\inn (a)}) - \sum_{l \in Q_{2}} 
\mc Hom (\mc V_{\out (l)}, \mc V_{\inn (l)})
- \sum_{i \in I} \mc End(\mc V_{i}), 
\]
and $\psi (\mc V)=\Eu^{\theta}(\Lambda_{Q}(\mc V))$ in 
$A^{\bullet}_{\mb T \times \C^{\ast}_{\theta} \times \C^{\ast}_{\hbar}}(\mc M)$.

\subsection{Computations of residues}
In \eqref{itcoho}, we have 
\[
\psi \left( \mc V \oplus \bigoplus_{i=1}^{j} \mc V^{(i)} \otimes e^{\hbar_{i}} \right) =
\psi(\mc V)  
\cdot
\prod_{i=1}^{j} 
\Eu^{\theta}\left( \mk N( \mc V^{(i)} \otimes e^{\hbar_{i}}, \mc V^{>i} 
) \right)
\cdot 
\prod_{i=1}^{j} 
\psi (\mc  V^{(i)} \otimes e^{\hbar_{i}}).  
\]
We divide the computations into three parts:
\begin{prop}
\label{compresidue}
We have the following.\\
(1) 
$
\displaystyle
\Res_{\hbar_{i} = \infty} 
\Eu^{\theta}\left( \mk N( \mc V^{(i)} \otimes e^{\hbar_{i}}, \mc V^{>i} ) \right)/
\Eu \left( \mk N( \mc V^{(i)} \otimes e^{\hbar_{i}}, \mc V^{>i} ) \right)
=
\bar{\beta}_{\infty} d_{i}\theta
$,
where
$d_{i}=|\mk I_{i}| /\beta_{0}$, and
\[
\displaystyle
\bar{\beta}_{\infty}=
\sum_{\substack{a \in Q_{1}\\ \out(a)=\infty}} 
\beta_{\inn(a)}  -
\sum_{\substack{a \in Q_{1}\\ \inn(a)=\infty}} 
\beta_{\out(a)}.
\]
(2) 
$
\int_{Fl_{H_{Q}( d_{i} \beta )}( \bar{\mc V}^{(i)}_{0}, \bar{\mk I}_{i})}
\psi (\mc  V^{(i)} \otimes e^{\hbar_{i}}) \Eu^{\theta}(\Theta( \mc F^{(i)}_{\bullet} )) 
= 
\frac{
(d_{i} \beta_{0} -1)!
}
{
\theta
}
 \int_{H_{Q}( d_{i} \beta )} 
\Eu^{\theta}(\Lambda_{Q^{\sharp}}(\mc V^{(i)} \oplus \mc V^{(i)}_{\infty'})).
$
\\
(3) 
$
\displaystyle
\Res_{\hbar = \infty} 
\Eu^{\theta}(\mk H(\mc F^{(i)}_{\bullet}, \mc F^{>i}_{\bullet} ))/
\Eu^{\theta}(\mk H(\mc F^{(i)}_{\bullet}, \mc F^{>i}_{\bullet} ))
 = 
- s( \mk I_{i}, \mbi{\mk I}_{>i}) \theta, 
$
where 
\[
s(\mk I_{i}, \mbi{\mk I}_{>i}) = 
\left| \lbrace (l, l') \in \mk I_{i} \times \mbi{\mk I}_{>i} \mid l < l' \rbrace \right|
-\left| \lbrace (l, l') \in \mk I_{i} \times \mbi{\mk I}_{>i} \mid l > l' \rbrace \right|.
\]
\end{prop}
\proof
(1)
We devide $\mathfrak{N}( \mc V^{\sharp}, \mc V^{\flat})$ into two parts: 
\begin{align*}
&
\sum_{(\heartsuit, \spadesuit)} 
\left(
\sum_{\substack{a \in Q_{1} \\ \inn(a), \out(a) \neq \infty}} 
\mc H om( \mc V^{\heartsuit}_{\out(a)}, \mc V^{\spadesuit}_{\inn(a)}) 
- 
\sum_{l \in Q_{2}} \mc H om( \mc V^{\heartsuit}_{\out(l)}, \mc V^{\spadesuit}_{\inn(\l)}) 
-
\sum_{i \in I} 
\mc Hom ( \mc V^{\heartsuit}_{i}, \mc V^{\spadesuit}_{i}) \right)\\
&+
\sum_{\substack{a \in Q_{1} \\ \out(a) = \infty}} 
\mc H om( \mc V^{\flat}_{\infty}, \mc V^{\sharp}_{\inn(a)}) 
+\sum_{\substack{a \in Q_{1} \\ \inn(a) = \infty}} 
\mc H om( \mc V^{\sharp}_{\out(a)}, \mc V^{\flat}_{\infty}). 
\end{align*}
The first line does not contribute to the residue, and dimension counting of the latter gives the assertion.
To prove (2) we remark that \eqref{ob4} implies
\begin{align*}
\Theta(\mc F^{(i)}_{\bullet})
&=
\mc Hom( \mc V^{(i)}_{\infty'}, \mc V^{(i)}_{0} ) - \mo_{M_{Q^{\sharp}}^{\zeta^{\sharp}} (d_{i} \beta)} 
+ \Theta(\bar{\mc F}^{(i)}_{\bullet}),
\end{align*}
and the right hand side does not include $e^{\hbar_{i}}$.
Since we have 
\begin{align*}
\Lambda_{Q^{\sharp}}(\mc V^{(i)} \oplus \mc V^{(i)}_{\infty'}) 
&=
\sum_{a \in Q_{1}} \mc Hom (\mc V^{(i)}_{\out (a)}, \mc V^{(i)}_{\inn (a)}) 
- \sum_{l \in Q_{2}} \mc Hom (\mc V^{(i)}_{\out (l)}, \mc V^{(i)}_{\inn (l)})
- \sum_{i \in I} \mc End(\mc V^{(i)}_{i})
\\
&+
\mc Hom( \mc V^{(i)}_{\infty'}, \mc V^{(i)}_{0}),
\end{align*}
the assertion follows from the integral
$\int_{Fl(\C^{d_{i}}, \bar{\mk I}_{(i)})} 
\Eu^{\theta}( T Fl(\C^{d_{i}}, \bar{\mk I}_{i}))=(d\beta_{0} -1)!$
over the compact manifold $Fl(\C^{d_{i}}, \bar{\mk I}_{i})$
which does not depend on $\theta$.
Finally (3) follows from direct computations.
\endproof

\subsection{Wall-crossing formula}
We compute \eqref{euler} by Proposition \ref{compresidue} when $\mk I=[\alpha_{0}]$.
We have
\begin{align}
&
\int_{\wt M^{\ell}(\alpha, [\alpha_{0}])} \tilde{\psi}
-\int_{M^{\zeta^{-}}(\alpha )} \psi 
\notag\\
&=
\sum_{\mk I^{\sharp} \in \mc D^{\ell}(\mk I)}
\frac{|\mk I^{\flat}|! (|\mk I_{\sharp}|-1)!}{\alpha_{0}!} 
\gamma_{d^{\sharp}}(\theta)
\cdot 
(s(\mk I^{\sharp}, \mk I^{\flat} ) - \bar{\beta}_{\infty} d^{\sharp} )
\cdot
\int_{\wt M^{ \min(\mk I^{\sharp})-1} (\alpha - d^{\sharp} \beta, \mk I^{\flat})} \tilde{\psi }.
\label{recursion}
\end{align}
where $\mk I^{\flat} = [\alpha_{0}] \setminus \mk I^{\sharp}$, $d^{\sharp} = |\mk I^{\sharp}| / \beta_{0}$, and
$\gamma_{d^{\sharp}} (\theta) =   
\int_{H_{Q} (d^{\sharp} \beta) }  
\Eu^{\theta}(\Lambda_{Q^{\sharp}}(\mc V^{\sharp} \oplus \mc V^{\sharp}_{\infty'}))$.

By Proposition \ref{compresidue}, we compute \eqref{main>0}, and we have
\begin{align}
&
\notag
\int_{M^{\zeta^{+}}(\alpha)} \Eu^{\theta}(\Lambda_{Q}(\mc V)) 
- \int_{M^{\zeta^{-}}(\alpha)} \Eu^{\theta}(\Lambda_{Q}(\mc V)) \\
\label{formula4}
&=
\sum_{j=1}^{\lfloor \alpha_{0} / \beta_{0} \rfloor} 
\sum_{\mbi{\mk I} \in \Dec{}_{\beta_{0}, j}^{\alpha_{0}}} 
\frac{ |\mk I_{\infty}|!  }{ \alpha_{0} ! }
\prod_{i=1}^{j} ( d_{i} \beta_{0} -1)! \gamma_{d_{i}} (\theta) ( s( \mk I_{i}, \mbi{\mk I}_{> i}) - \bar{\beta}_{\infty} d_{i} ) 
\int_{ M^{\zeta^{-}}(\alpha - |\mbi d_{\mbi{\mk I}}| \beta )} \Eu^{\theta}(\Lambda_{Q}(\mc V)),
\end{align}
where $d_{i} = |\mk I_{i}|/\beta_{0}$ and 
$\mbi{\mk I}_{> i} = \bigsqcup_{k >i} \mk I_{k} \sqcup \mk I_{\infty}$ 
for $\mbi{\mk I}=( \mk I_{1}, \ldots, \mk I_{j}) \in \text{Dec}_{\beta_{0}, j}^{\alpha_{0}}$, and
$\gamma_{d_{i}} (\theta) =   
\int_{H_{Q} (d_{i} \beta) }  
\Eu^{\theta}(\Lambda_{Q^{\sharp}}(\mc V^{\sharp} \oplus \mc V^{\sharp}_{\infty'}))$.

We put $\Dec (\alpha_{0})= \bigsqcup_{j=1}^{\lfloor \alpha_{0}/\beta_{0} \rfloor} 
\Dec{}_{\beta_{0}, j}^{\alpha_{0}}$ and, we get the following theorem.
\begin{thm}
\label{thm:adjoint}
\begin{align*}
&
\int_{M^{\zeta^{+}}(\alpha)} \Eu^{\theta}(\Lambda_{Q}(\mc V)) 
- \int_{M^{\zeta^{-}}(\alpha)} \Eu^{\theta}(\Lambda_{Q}(\mc V)) \\
\label{formula4}
&=
\sum_{k=0}^{\lfloor \alpha_{0}/\beta_{0} \rfloor} 
\sum_{\substack{\mbi{\mk I} \in \Dec(\alpha_{0}) \\ |\mbi{\mk I}|=k}} 
\frac{ |\mk I_{\infty}|!  }{ \alpha_{0} ! }
\prod_{i=1}^{j} ( d_{i} \beta_{0} -1)! \gamma_{d_{i}} (\theta) ( s( \mk I_{i}, \mbi{\mk I}_{> i}) - \bar{\beta}_{\infty} d_{i} ) 
\int_{ M^{\zeta^{-}}(\alpha -k \beta)} \Eu^{\theta}(\Lambda_{Q}(\mc V)).
\end{align*}
\end{thm}

\begin{cor}
\label{cor:indep}
When $Q$ is as in Example \ref{exam:quivervar}, \ref{exam:blup}, then the coefficients of wall-crossing term does not depend
on $r$.
\end{cor}

\begin{NB}
For $\mbi{\mk I} =(\mk I_{1}, \ldots, \mk I_{j})$, we put 
\[
a_{\mbi{\mk I}}(\bar{\beta}_{\infty})=\prod_{i=1}^{j} 
( s(\mk I_{i}, \mbi{\mk I}_{> i}) - \bar{\beta}_{\infty} d_{i} ) (|\mk I_{i}|-1)! \gamma_{d_{i}}(\theta).
\]
We also write by $\Dec (i, \alpha_{0})$ the set of elements 
$\mbi{\mk I}
\in \bigsqcup_{j=1}^{\lfloor \alpha_{0}/\beta_{0} \rfloor} 
\Dec{}_{\beta_{0}, j}^{\alpha_{0}}$ such that $|\mk I_{\infty}|=i$, and put 
\[
A_{\alpha_{0}} ( \bar{\beta}_{\infty})= \sum_{\mbi{\mk I} \in \Dec(0, \alpha_{0})} a_{\mbi{\mk I}}( \bar{\beta}_{\infty}).
\]
\begin{ques}
For any fixed $i$, can we show that 
$\sum_{\mbi{\mk I} \in \Dec(i, \alpha_{0})} a_{\mbi{\mk I}} (\bar{\beta}_{\infty})= 
\begin{pmatrix} \alpha_{0} \\ i \end{pmatrix} 
A_{\alpha_{0}-i}( \bar{\beta}_{\infty})$?
\end{ques}

If this is correct, substituting these into \eqref{formula4}, we get the following theorem.
\begin{thm}
\label{thm:adjoint2}
\begin{align*}
&
\int_{M^{\zeta^{+}}(\alpha)} \Eu^{\theta}(\Lambda_{Q}(\mc V)) 
- \int_{M^{\zeta^{-}}(\alpha)} \Eu^{\theta}(\Lambda_{Q}(\mc V)) \\
\label{formula4}
&=
\sum_{i=0}^{\alpha_{0}} \frac{ i!  }{ \alpha_{0}! } \sum_{\mbi{\mk I} \in \Dec(i, \alpha_{0})} 
a_{\mbi{\mk I}} (\bar{\beta}_{\infty} )
\int_{ M^{\zeta^{-}}(\alpha - |\mbi d_{\mbi{\mk I}}| \beta )} \Eu^{\theta}(\Lambda_{Q}(\mc V))
\\
&=
\sum_{k=0}^{\lfloor \alpha_{0}/\beta_{0} \rfloor} \frac{A_{k\beta_{0}}( \bar{\beta}_{\infty})}{(k \beta_{0})!} 
\int_{ M^{\zeta^{-}}(\alpha -k \beta)} \Eu^{\theta}(\Lambda_{Q}(\mc V)).
\end{align*}
\end{thm}

We can explicitly write $A_{\alpha_{0}} (\bar{\beta}_{\infty})$ as 
\begin{align*}
A_{\alpha_{0}} (\bar{\beta}_{\infty}) 
&= 
\sum_{\substack{\mk I_{1} \sqcup \cdots \sqcup \mk I_{j} = [\alpha_{0}] \\
\min(\mk I_{1}) > \cdots > \min(\mk I_{j})\\
|\mk I_{i}| = d_{i} \beta_{0} }}
\prod_{i=1}^{j} (s( \mk I_{i}, \mbi{\mk I}_{>i} ) - \bar{\beta}_{\infty} d_{i} )(d_{i} -1)! \gamma_{d_{i}} (\theta)
\end{align*}
For a fixed $d=(d_{1}, \ldots, d_{j})$, we check coefficients of 
$\prod_{i=1}^{j} (d_{i} -1)! \gamma_{d_{i}} (\theta)$ vanish when $\bar{\beta}_{\infty}=0$.

Namely we check 
\begin{align*}
\sum_{\substack{\mk I_{1} \sqcup \cdots \sqcup \mk I_{j} = [\alpha_{0}] \\
\min(\mk I_{1}) > \cdots > \min(\mk I_{j})\\
|\mk I_{i}| = \beta_{0} d_{i}  }}
\prod_{i=1}^{j} s( \mk I_{i}, \mbi{\mk I}_{>i} )
=0.
\end{align*}

\begin{verbatim}
EWop:=(y,m)->[op(y),m]:

ListSetPtns:=proc(n,k)
local east,west,i,out ; options remember;
#Lists all partitions of the set 1,..,n into k classes #Output array[i] is the class to which i belongs.
if n=1 then
if k<>1 then RETURN([]) else RETURN([[1]]) fi: else
east:=ListSetPtns(n-1,k-1): west:=ListSetPtns(n-1,k): out:=map(EWop,east,k);
for i from 1 to k do out:=[op(out),op(map(EWop,west,i))] od; 
RETURN(out);
fi: end:






dimv := proc(L, j) 
local i, M; 
M := [seq(0, i = 1 .. j )]; 
for i from 1 to nops(L) do 
if L[i] <= j then
M := applyop(x -> x + 1, L[i], M); 
fi;
end do; 
return M; 
end;




s := proc(L, l) 
local i, j, x, y; 
x := 0; y := 0; 
for i from 1 to nops(L) - 1 do for j from i + 1 to nops(L) do 
if L[i] = l then if l < L[j]  then x := x + 1; end if; end if; 
if L[j] = l then if l < L[i] then y := y + 1; end if; end if; 
end do; end do; 
RETURN(x-y); 
end;


check:=proc(d, n)
local  i,j, f, k, l, m, L, M;
# d=(d_{1}, \ldots, d_{j})
j:=nops(d);
f:=0;
for L in ListSetPtns(n, j+1) do 
for k from 1 to j+1 do
M:=[];
for i from 1 to n do
if L[i] < k then
M:=[op(M), L[i]];
fi;
if L[i] > k then
M:=[op(M), L[i]-1];
fi;
if L[i] = k then
M:=[op(M), j+1];
fi;
od;
if dimv(M,j)= d then
m:=nops(M);
M:=[seq(M[m-i+1], i=1..m)];
f:=f+product(s(M, l), l=1..j);
fi;
od;
od;
for L in ListSetPtns(n, j) do 
if dimv(L,j)= d then
L:=[seq(L[n-i+1], i=1..n)];
f:=f+product(s(L, l), l=1..j);
fi;
od;
RETURN(factor(f));
end;




\end{verbatim}

\end{NB}

\subsection{When $Q$ has only one arrow from $\infty$}
Motivated by $Q^{\sharp}$-representations studied in \S \ref{subsec:moduli}, we consider the following situation.
Here we assume that we have only one arrow connecting with $\infty$ in our quiver $Q$, and it starts form $\infty$ and ends at $0 \in I$.
We also assume $\out(l), \inn(l) \in I$ for any $l \in Q_{2}$ as so far.

Then in \eqref{recursion} we have $\wt M^{\ell}(\alpha) = \emptyset$ for $\ell>0$ since 
we can always take a sub-representation $S=\bigoplus_{i \in I} V_{i}$ with $S_{\infty}=0$ in Definition \ref{ellstab}.
This also implies that in the right hand side of \eqref{recursion} it is enough to consider
$\mk I^{\sharp} \subset [r]$ containing $1$.

We put $a_{\alpha}=\int_{M^{\zeta^{-}}( \alpha )} \Eu^{\theta}(\Lambda_{Q}(\mc V))$. 
Then \eqref{recursion} is re-written as 
\begin{align*}
a_{\alpha}
&=
\sum_{\substack{1 \in \mk I^{\sharp} \subset [\alpha_{0}] \\ |\mk I^{\sharp}|/\beta_{0} \in  \Z_{>0} } }
\frac{|\mk I^{\flat}|! (|\mk I^{\sharp} | - 1)!}{\alpha_{0}!}
\gamma_{d^{\sharp}}(\theta)
\cdot 
(- s(\mk I^{\sharp}, \mk I^{\flat}) + \beta_{0} d^{\sharp} )
\cdot
a_{\alpha - d^{\sharp} \beta}
\end{align*}
where $d^{\sharp} = |\mk I^{\sharp}|/\beta_{0}$ and $\mk I^{\flat} = [\alpha_{0}] \setminus \mk I^{\sharp}$.
In particular, when $\alpha=d\beta$ we have $a_{d \beta} = \gamma_{d}(\theta)$ since $Q=Q^{\sharp}$ 
and $M_{Q}^{\zeta}(d \beta) \cong H_{Q}(d \beta)$. 
Hence we have
\begin{align*}
\gamma_{d}( \theta)
&=
\sum_{\substack{1 \in \mk I^{\sharp} \subset [d\beta_{0}] \\ |\mk I^{\sharp}|/\beta_{0} \in  \Z_{>0} } }
\frac{|\mk I^{\flat}|! (|\mk I^{\sharp} | - 1)!}{(d\beta_{0})!}
(- s(\mk I^{\sharp}, \mk I^{\flat}) + \beta_{0} d^{\sharp} )
\cdot
\gamma_{d^{\sharp}}(\theta)
\cdot 
\gamma_{d- d^{\sharp}}(\theta).
\end{align*}

\begin{NB}
When we consider the $A_{1}$-handsaw quiver with the framing $(r_{0}, r_{1})=(1,0)$, we have
$\gamma_{d}(\theta)= (\theta/\e+1)_{d}/ d!$.
We have
\begin{align*}
0
&=
\sum_{\substack{1 \in \mk I^{\sharp} \subsetneq [d]  } }
\frac{|\mk I^{\flat}|! (|\mk I^{\sharp} | - 1)!}{(d\beta_{0})!}
(- s(\mk I^{\sharp}, \mk I^{\flat}) + \beta_{0} d^{\sharp} )
\cdot
\frac{(\theta/\e+1)_{d} }{d!} \cdot 
\frac{(\theta/\e+1)_{d-d^{\sharp}} }{(d-d^{\sharp})!}.
\end{align*}

We check this by the following maple computations :

\begin{verbatim}
EWop:=(y,m)->[op(y),m]:

Lsp:=proc(n,k)
local east,west,i,out ; options remember;
#Lists all partitions of the set 1,..,n into k=2 classes C_{1}=\mk I_{\sharp} and C_{2}=\mk I_{\flat} 
#Output array[i] is the class to which i belongs. i \in C_{array[i]}
#Always array[1]=1, i.e. 1 \in C_{1} 
if n=1 then
if k<>1 then RETURN([]) else RETURN([[1]]) fi: else
east:=Lsp(n-1,k-1): west:=Lsp(n-1,k): out:=map(EWop,east,k);
for i from 1 to k do out:=[op(out),op(map(EWop,west,i))] od; RETURN(out);
fi: end:


s:=proc(L)
local i, j, x;
x:=0;
for i from 1 to nops(L)-1 do
for j from i+1 to nops(L) do
if L[i] <> L[j] then
if L[i] =1 then x:=x+1; else x:= x-1; fi; else fi;
od; od;
RETURN(x)
end:

poch:=proc(x, k)
local f, i;    
f:=1; for i from 0 to k-1
do f:=f*(x+i) od;
RETURN(f);
end:

check:=proc(n)
local i, L, d, f;
f:=0;
for L in Lsp(n,2) do 
d:=0;
for i from 1 to nops(L) do 
if L[i]=1 then d:=d+1;
else fi;
od;
f:=f+factorial(n-d)*factorial(d-1)/factorial(n)*( -s(L)+d)*poch(theta+1,d)/poch(1, d)*poch(theta+1,n-d)/poch(1, n-d)
od;
RETURN(factor(f))
end;


\end{verbatim}

\end{NB}


\begin{thebibliography}{NY}
\bibitem{AHKOSSY}
H.~Awata, K.~Hasegawa, H.~Kanno, R.~Ohkawa, S.~Shakirov, J.~Shiraishi and Y.~Yamada
{\it Non-stationary difference equation, affine Laumon space
and quantization of discrete Painlev\'e equation}, arXiv:2211.16772
%


\bibitem{AGT}L. Alday, D. Gaiotto and Y. Tachikawa, {\it Liouville correlation functions from four- dimensional gauge theories}, Lett. Math. Phys. {\bf 91} 167--197 (2010)

\bibitem{BE} A. Braverman and P. Etingof, {\it Instanton counting via affine Lie algebras II: from Whittaker vectors to the Seiberg-Witten prepotential}, In: Bernstein J., Hinich V., Melnikov A., (eds), Studies in Lie theory, Progr. Math. 243,  61--78, Birkhauser, Boston (2006)

\bibitem{BF} K. Behrend and B. Fantechi, {\it The intrinsic normal cone}, Invent. Math.{\bf 128}, 45--88 (1997)

\bibitem{BS1}M. Bershtein and A. Shchechkin, {\it Bilinear equations on Painlev\'e $\tau$ functions from CFT}, Comm. Math. Phys. 339 (2015), no. 3, 1021--1061. 
 
 
 
\bibitem{C} Crawley-Boevey, {\it Geometry of the moment map for representations of quivers}, Compositio Math.{\bf 126}, 257--293 (2001)


\bibitem{F}W. Fulton, {\it Intersection Theory}, Springer-Verlag, Berlin (1984)


\bibitem{GNY} L. G\"{o}ttsche, H. Nakajima and K. Yoshioka, {\it Donaldson = Seiberg-Witten from Mochizuki's formula and Instanton Counting}, Publ. RIMS Kyoto Univ.{\bf 47}, 307--359 (2011)

\bibitem{GP} T.Graber and R. Pandharipande, {\it Localization of virtual classes}, Invent. Math. 135  (1999), 487--518.

%

\bibitem{Hasegawa}
K.~Hasegawa,
{\it Quantizing the B\"acklund transformations of Painlev\'e equations and the quantum discrete Painlev\'e VI equation},
{\it Adv. Stud. Pure Math.} \textbf{61}
{\it ``Exploring New Structures and Natural Constructions in Mathematical Physics"}, K. Hasegawa, et. al, eds.
(2011), 275-288

\bibitem{HasegawaLax}
K.~Hasegawa,
{\it Quantizing the Discrete Painleve VI Equation : The Lax Formalism},
Lett. Math. Phys. \textbf{103} (2013) 865-879



\bibitem{HLNR1} M. Halln\"as, E. Langmann, M. Noumi and H. Rosengren,
{\it From Kajihara’s transformation formula to deformed Macdonald-Ruijsenaars and Noumi-Sano operators}, 
Selecta Math. (N.S.) 28 (2022), no. 2, Paper no. 24, 36 pp.

\bibitem{HLNR2} M. Halln\"as, E. Langmann, M. Noumi and H. Rosengren,
{\it Higher order deformed elliptic Ruijsenaars operators}, 
Comm. Math. Phys. 392 (2022), no. 2, 659--689. 
 

%


\bibitem{IMO} Y. Ito, K. Maruyoshi, T. Okuda, {\it Scheme dependence of instanton counting in ALE spaces}, J. High Energy Phys. 2013, no. 5, 045, front matter+16 pp.

\bibitem{K} Y. Kajihara, 
{\it Euler transformation formula for multiple basic hypergeometric series of
type $A$ and some applications}, Adv. in Math. 187 (2004), 53--97.

\bibitem{KN} Y. Kajihara, M. Noumi, 
{\it Multiple elliptic hypergeometric series. 
An approach from the Cauchy determinant}, Indag. Mathem., N. S., 14 (3, 4), 395--421.

\bibitem{Kr} A. Kresch, {\it Cycle groups for Artin stacks}, Invent. Math.{\bf 138}, 495--536 (1999) 

%

\bibitem{LSW} R. Langer, M.J. Schlosser and S.O. Warnaar, 
{\it Theta functions, elliptic hypergeometric series, and Kawanaka's Macdonald polynomial conjecture}, 
SIGMA 5 (2009), Paper 055.

\bibitem{M}T. Mochizuki, {\it Donaldson Type Invariants for Algebraic Surfaces: Transition of Moduli Stacks}, Lecture Notes in Math. 1972, Springer, Berlin, 2009.

%
\bibitem{N1}H. Nakajima, {\it Instantons on ALE spaces, quiver varieties, and Kac-Moody algebras}, Duke Math. J. 76 (1994), no. 2, 365--416.
%
%
%
 
\bibitem{Nek} N. Nekrasov, {\it Seiberg-Witten prepotential from instanton counting}, Adv. Theor. Math. Phys.{\bf 7}, no. 5, 831--864 (2003)

\bibitem{NO} N. Nekrasov and A. Okounkov, {\it Seiberg-Witten theory and random partitions}, In: Etingof P., Retakh V. S., Singer, I.M., (eds), The unity of mathematics, Progr. Math. 244, 525--596, Birkhauser Boston, Boston, MA, (2006)
%

\bibitem{NY1} H. Nakajima and K. Yoshioka,
{\it Instanton counting on blowup. I. $4$-dimensional pure gauge theory}, Invent. Math. 162 (2005), no. 2, 313--355

\bibitem{NY2} H. Nakajima and K. Yoshioka, {\it Perverse coherent sheaves on blowup. III. Blow-up formula from wall-crossing}, Kyoto Journal of Mathematics, Vol. 51, No. 2 (2011), 263--335

\bibitem{O1} R. Ohkawa, {\it Wall-crossing between stable and co-stable ADHM data}, Lett. Math. Phys.  108 (2018), no. 6 1485--1523

\bibitem{O2}R. Ohkawa, {\it Functional Equations of Nekrasov Functions Proposed by Ito, Maruyoshi, and Okuda},
Moscow Math. J. 20 (2020), no. 3, 531--573

\bibitem{OY}R. Ohkawa, Y. Yoshida,
{\it Wall-crossing for vortex partition function and handsaw quiver varierty}, arXiv:2208.00435

\bibitem{O3} R. Ohkawa, {\it Residue formula for flag manifold of type $A$ from wall-crossing}, in preparation

\bibitem{R} M. Reineke, {\it Poisson automorphisms and quiver moduli}, J. Inst. Math. Jussieu 9 (2010), no. 3, 653--667

\bibitem{Shakirov:2021krl}
S.~Shakirov,
{\it Non-stationary difference equation for $q$-Virasoro conformal blocks},
arXiv:2111.07939

\bibitem{S}A. Shchechkin, {\it Blowup relations on $\C^{2}/\Z_{2}$ from Nakajima-Yoshioka blowup relations},
arXiv:2006.08582

\bibitem{Shi}
J.~Shiraishi,
{\it Affine screening operators, affine Laumon spaces and conjectures
  concerning non-stationary Ruijsenaars functions},
  J.Integrable.Syst. {\bf 4} (2019) xyz010



%
\end{thebibliography}
\end{document}